\documentclass[10pt]{amsart}
\usepackage{amsfonts}
\usepackage{amsthm}
\usepackage{amsmath}
\usepackage{graphicx}
\usepackage{amscd,amssymb,amsthm}
\usepackage{enumerate}
\usepackage{mathrsfs}
\usepackage{epstopdf}
\usepackage{mathptmx}

\newcounter{minutes}
\divide\time by 60
\newcounter{hours}
\multiply\time by 60 \addtocounter{minutes}{-\time}

\title[Cross-product of Bessel functions]{Cross-product of Bessel functions: monotonicity patterns and functional inequalities$^{\bigstar}$}
\thanks{$^{\bigstar}$The work of \'A. Baricz was supported by the J\'anos Bolyai Research Scholarship of
the Hungarian Academy of Sciences. The second author is on leave from  IIT Madras. The research of S. Singh was supported by the
fellowship of the University Grants Commission, India. The authors are grateful to Prof. Tibor K. Pog\'any
for a fruitful discussion on the zeros of the cross-product of Bessel functions
during his visit to Indian Statistical Institute, Chennai Centre in February 2015.}

\author[\'A. Baricz]{\'Arp\'ad Baricz}
\address{Institute of Applied Mathematics, \'Obuda University, 1034 Budapest, Hungary}
\address{Department of Economics,  Babe\c{s}-Bolyai University, Cluj-Napoca 400591, Romania}
\email{bariczocsi@yahoo.com}

\author[S. Ponnusamy]{Saminathan Ponnusamy}
\address{Indian Statistical Institute, Chennai Centre, Society for Electronic Transactions and Security,
MGR Knowledge City, CIT Campus, Taramani, Chennai 600113, India}
\email{samy@iitm.ac.in, samy@isichennai.res.in}

\author[S. Singh]{Sanjeev Singh}
\address{Department of Mathematics, Indian Institute of Technology Madras, Chennai 600036, India}
\email{sanjeevsinghiitm@gmail.com}

\newtheorem{theorem}{Theorem}

\newtheorem{lemma}{Lemma}
\newtheorem{corollary}{Corollary}

\begin{document}

\def\thefootnote{}
\footnotetext{ \texttt{File:~\jobname .tex,
          printed: \number\year-0\number\month-0\number\day,
          \thehours.\ifnum\theminutes<10{0}\fi\theminutes}
} \makeatletter\def\thefootnote{\@arabic\c@footnote}\makeatother



\maketitle


\begin{abstract}
In this paper we study the Dini functions and the cross-product of Bessel functions. Moreover, we are interested on the monotonicity patterns for the cross-product of Bessel and modified Bessel functions. In addition, we deduce Redheffer-type inequalities, and the interlacing property of the zeros of Dini functions and the cross-product of Bessel and modified Bessel functions. Bounds for logarithmic derivatives of these functions are also derived. The key tools in our proofs are some recently developed infinite product representations for Dini functions and cross-product of Bessel functions.
\end{abstract}

\section{\bf Introduction and preliminaries}
\setcounter{equation}{0}

Bessel and modified Bessel functions of the first kind play an important role in the theory of special functions because they are useful in many problems of applied mathematics. These functions have been studied by many researchers, and their study goes back to famous scientists like Bessel, Euler, Fourier, and others. Motivated by their appearance as eigenvalues in the clamped plate problem for the ball, Ashbaugh and Benguria have conjectured
that the positive zeros of the cross-product of Bessel and modified Bessel functions of first kind, defined by
$$
W_{\nu}(z)=J_{\nu}(z)I'_{\nu}(z)-J'_{\nu}(z)I_{\nu}(z)=J_{\nu+1}(z)I_{\nu}(z)+J_{\nu}(z)I_{\nu+1}(z),
$$
where $J_{\nu}$ and $I_{\nu}$ stand for the Bessel and modified Bessel functions of the first kind, increase with $\nu$ on $\left[-\frac{1}{2},\infty\right).$ Lorch \cite{lorch} verified this conjecture and presented some other properties of the zeros of the above cross-product of Bessel and modified Bessel functions. See also the paper \cite{ben} of Ashbaugh and Benguria for more details. Recently, the authors of \cite{abp} pointed out that actually the above monotonicity property is valid on $(-1,\infty)$ and proved that for $\nu>-1$ and $z\in \mathbb{C}$ the power series representation
\begin{equation}\label{sum}
W_{\nu}(z)=2\sum_{n\ge 0}\frac{(-1)^n(\frac{z}{2})^{2\nu+4n+1}}{n!\Gamma(\nu+n+1)\Gamma(\nu+2n+2)}
\end{equation}
and the infinite product representation
\begin{equation}\label{product}
\mathscr{W}_{\nu}(z)=2^{2\nu}z^{-2\nu-1}\Gamma(\nu+1)\Gamma(\nu+2)W_{\nu}(z)=\prod_{n\geq 1}\left(1-\frac{z^4}{\gamma^4_{\nu,n}}\right)
\end{equation}
are valid, see \cite{abp} for more details. In this paper we would like to continue the study of the properties of the cross-product of Bessel and modified Bessel functions of the first kind by showing a series of new results. We also consider another special combination of Bessel functions, namely, the so called Dini functions $d_{\nu}:\mathbb{C}\rightarrow \mathbb{C},$ defined by
$$d_{\nu}(z)=(1-\nu)J_{\nu}(z)+zJ'_{\nu}(z)=J_{\nu}(z)-zJ_{\nu+1}(z),$$
and the modified Dini function $\xi_{\nu}:\mathbb{C}\to\mathbb{C},$ defined by
$$\xi_{\nu}(z)=\mathrm{i}^{-\nu}d_{\nu}(\mathrm{i}z)=(1-\nu)I_{\nu}(z)+zI'_{\nu}(z)
=I_{\nu}(z)+zI_{\nu+1}(z).$$
For $\nu>-1$ and $z\in \mathbb{C}$, the Weierstrassian factorization of Dini functions is \cite{bpogsz}
\begin{equation}\label{dfprod}
\mathcal{D}_{\nu}(z)=2^{\nu}\Gamma(\nu+1)z^{-\nu}d_{\nu}(z)=\prod_{n\geq 1}\left(1-\frac{z^2}{\alpha_{\nu,n}^2}\right)
\end{equation}
and the Weierstrassian factorization of modified Dini function is \cite{bps2}
\begin{equation}\label{mdfprod}
\lambda_{\nu}(x)=2^{\nu}\Gamma(\nu+1)x^{-\nu}\xi_{\nu}(x)=\prod_{n\geq 1}\left(1+\frac{x^2}{\alpha_{\nu,n}^2}\right).
\end{equation}
The paper is organized as follows: the next section contains the main results on the cross-product of Bessel functions, Dini functions and their zeros. Section 3 is devoted for the proofs of the main results. In our proofs we use a series of methods: Mittag-Leffler expansions, Laguerre separation theorem, Laguerre inequality for entire functions, differential equation for the Dini function, monotone form of L'Hospital's rule, and representations of logarithmic derivatives of Dini functions and cross-product of Bessel functions via the spectral zeta functions of the zeros of the above functions.

\section{\bf Main results}\label{MR}
\setcounter{equation}{0}

\subsection{Monotonicity properties} Our first set of results are some monotonicity and concavity properties of cross-product of Bessel and modified Bessel functions of first kind.

\begin{theorem}\label{theorem1}
Let $\nu>-1$ and define $S=S_1\cup S_2$, where
$S_1=\bigcup_{n\geq 1}[-\gamma_{\nu,2n},-\gamma_{\nu,2n-1}]$, $S_2=\bigcup_{n\geq 1}[\gamma_{\nu,2n-1},\gamma_{\nu,2n}]$
and $\gamma_{\nu,n}$ denotes the $n$th positive zero of the function $W_{\nu}$.
Then the following assertions hold true:
\begin{enumerate}
\item[\bf a.]The function $x\mapsto\mathscr{W}_{\nu}(x)$ is negative on $S$ and it is strictly positive on $\mathbb{R}\setminus S;$
\item[\bf b.] The function $x\mapsto\mathscr{W}_{\nu}(x)$ is strictly increasing on $(-\gamma_{\nu,1},0]$ and strictly decreasing on $[0,\gamma_{\nu,1});$
\item[\bf c.] The function $x\mapsto\mathscr{W}_{\nu}(x)$ is strictly log-concave on $\mathbb{R}\setminus S$ and
strictly geometrically concave on $(0,\infty)\setminus S_2;$
\item[\bf d.] The function $x\mapsto W_{\nu}(x)$ is strictly log-concave on $(0,\infty)\setminus S_2$ for all $\nu\geq-\frac{1}{2};$
\item[\bf e.] The function $\nu\mapsto \mathscr{W}_{\nu}(x)$ is increasing on $(-1,\infty)$ for all $x\in (-\gamma_{\nu,1},\gamma_{\nu,1})$ and the function $\nu\mapsto x\mathscr{W}'_{\nu}(x)/\mathscr{W}_{\nu}(x)$ is increasing on $(-1,\infty)$ for all $x\in \mathbb{R};$
\item[\bf f.] The function $x\mapsto\left(-\log\mathscr{W}_{\nu}(\sqrt[4]{x})\right)'$ and the
function $x\mapsto1/\mathscr{W}_{\nu}(\sqrt[4]{x})$ are absolutely monotonic on $(0,\gamma^4_{\nu,1})$ for all $\nu>-1.$
\end{enumerate}
\end{theorem}
\subsection{Interlacing of positive real zeros of Bessel and related functions}
Let us recall Dixon's theorem \cite[p.480]{watson} which states that, when $\nu>-1$ and $a,b,c,d$ are
constants such that $ad\neq bc$, then the positive zeros of $x\mapsto aJ_{\nu}(x)+bxJ'_{\nu}(x)$ are interlaced with those of
$x\mapsto cJ_{\nu}(x)+dxJ'_{\nu}(x)$. Therefore if we choose $a=1-\nu$, $b=c=1$ and $d=0$ then for $\nu>-1$ we have,
\begin{equation}\label{IP1}
j_{\nu,n-1}<\alpha_{\nu,n}<j_{\nu,n}, ~~\mbox{where}~~ n\in \mathbb{N},
\end{equation}
with the convention that $j_{\nu,0}=0$. Here $j_{\nu,n}$ stands for the $n$th positive zero of the Bessel function $J_{\nu}.$ In \cite{abp}, among other things, the following interlacing inequality has been proved for $\nu>-1$,
\begin{equation}\label{IP2}
j_{\nu,n}<\gamma_{\nu,n}<j_{\nu,n+1}, ~~\mbox{where}~~ n\in \mathbb{N}.
\end{equation}
Taking into account the above two interlacing inequalities, it is natural to ask whether the zeros of Dini functions and of the cross-product of Bessel functions satisfy some interlacing property. The next theorem will answer this question.

\begin{theorem}\label{theorem2}
For $\nu>-1$, the zeros of Dini functions and of the cross-product of Bessel functions are interlacing, that is, they satisfy the following interlacing inequality
\begin{equation}\label{IP3}
\alpha_{\nu,n}<\gamma_{\nu,n}<\alpha_{\nu,n+1}, ~~\mbox{where}~~ n\in \mathbb{N}.
\end{equation}
\end{theorem}

Thus, combining the inequalities \eqref{IP1}, \eqref{IP2} and \eqref{IP3} we have the following:
$$\alpha_{\nu,n}<j_{\nu,n}<\gamma_{\nu,n}<\alpha_{\nu,n+1}<j_{\nu,n+1}, ~~\mbox{where}~~ n\in \mathbb{N}.$$
As an immediate consequence of the above interlacing properties we have the following upper and lower bounds for the cross-product of Bessel functions and consequently we can get bounds for ratio of modified Bessel and Bessel functions of first kind.
\begin{corollary}\label{corol2.1}
If $\nu>-1$, then the following inequalities hold:
\begin{equation}\label{corol2.1_ine1}
\mathcal{D}_{\nu}(x)\lambda_{\nu}(x)<\mathscr{W}_{\nu}(x)<\frac{\alpha^4_{\nu,1}}{\alpha^4_{\nu,1}-x^4}\mathcal{D}_{\nu}(x)\lambda_{\nu}(x)\ \ \ \ \mbox{for}\ \ |x|<\alpha_{\nu,1}
\end{equation}
\begin{equation}\label{corol2.1_ine2}
\mathcal{J}_{\nu}(x)\mathcal{I}_{\nu}(x)<\mathscr{W}_{\nu}(x)<\frac{j^4_{\nu,1}}{j^4_{\nu,1}-x^4}\mathcal{J}_{\nu}(x)\mathcal{I}_{\nu}(x) \ \ \ \ \mbox{for}\ \ |x|<j_{\nu,1}
\end{equation}
and
\begin{equation}\label{corol2.1_ine3}
e^{\frac{x^2}{2(\nu+1)}}<\frac{I_{\nu}(x)}{J_{\nu}(x)}<\left(\frac{j^2_{\nu,1}+x^2}{j^2_{\nu,1}-x^2}\right)^{\frac{j^2_{\nu,1}}{4(\nu+1)}}\ \ \ \ \mbox{for}\ \ 0<x<j_{\nu,1}.
\end{equation}
The reverse inequalities in \eqref{corol2.1_ine3} holds for $-j_{\nu,1}<x<0$.
\end{corollary}
In view of the inequality $\alpha_{\nu,n}<j_{\nu,n}$ where $n\in \mathbb{N},$ we observe that the left-hand side inequality of \eqref{corol2.1_ine2} is better than the left-hand side inequality of \eqref{corol2.1_ine1} while the right-hand side inequality of \eqref{corol2.1_ine1} is better than the right-hand side inequality of \eqref{corol2.1_ine2}.

Moreover, the next interlacing properties are also valid.

\begin{theorem}\label{theorem3}
If $\nu>-1$, then the following interlacing properties are valid:
\begin{enumerate}
\item[\bf a.] The zeros of the function $z\mapsto \mathcal{D}'_{\nu}(z)$ are interlaced with those of the function $z\mapsto \mathcal{D}_{\nu}(z)$.
\item[\bf b.] The zeros the function $z\mapsto \mathscr{W}'_{\nu}(\sqrt{z})$ are interlaced with those of the function $z\mapsto\mathscr{W}_{\nu}(\sqrt{z})$.
\item[\bf c.] For $\nu>0$, the zeros of the function $z\mapsto d'_{\nu}(z)$ are interlaced with those of the function $z\mapsto d_{\nu}(z)$.
\end{enumerate}
\end{theorem}

Now, we present an identity for zeros of Dini functions and zeros of cross-product of Bessel functions which is analogous to the identity of Calogero for the zeros of Bessel functions of the first kind, see \cite{bmps} for more details.

\begin{theorem}\label{theorem4}
Let $\nu>-1$, $k\in \mathbb{N}$ and $\gamma'_{\nu,k}$ denote the $k$th positive zero of the function $W'_{\nu}$. Then the following identities are valid:
\begin{equation}\sum_{n\ge 1, n\neq k}\frac{1}{\alpha^2_{\nu,n}-\alpha^2_{\nu,k}}=\frac{1}{4\alpha^2_{\nu,k}}
\left[2\nu+1-\frac{\alpha^2_{\nu,k}+2\nu-1} {\alpha^2_{\nu,k}-2\nu+1}\right], \label{id1} \end{equation}
\begin{equation}\sum_{n\ge 1, n\neq k}\frac{1}{\alpha^4_{\nu,n}-\alpha^4_{\nu,k}}=\frac{1}{8\alpha^4_{\nu,k}}
\left[2\nu+3-\frac{\alpha^2_{\nu,k}+2\nu-1}{\alpha^2_{\nu,k}-2\nu+1}\right]-\frac{1}{2\alpha^2_{\nu,k}}\sum_{n\ge 1}\frac{1}{\alpha^2_{\nu,n}+\alpha^2_{\nu,k}}, \label{id2} \end{equation}
and
\begin{equation}\sum_{n\ge 1, n\neq k}\frac{1}{\gamma^4_{\nu,n}-\gamma^4_{\nu,k}}=\frac{1}{8\gamma^4_{\nu,k}}
\left[2\nu+5+\sum_{n\ge 1}\frac{4\gamma^4_{\nu,k}}{\gamma'^4_{\nu,n}-\gamma^4_{\nu,k}}\right]. \label{id3}
\end{equation}
\end{theorem}

\subsection{\bf Rayleigh functions}
Before we state our next result, let us define the Rayleigh functions (or spectral zeta functions) for the zeros of Dini function and for the zeros of cross-product of Bessel and modified Bessel functions by
\begin{equation}\label{rf1}
\eta_{2m}(\nu)=\sum_{n\geq 1}\frac{1}{\alpha^{2m}_{\nu,n}}
\end{equation}
and
\begin{equation}\label{rf2}
\zeta_{4m}(\nu)=\sum_{n\geq 1}\frac{1}{\gamma^{4m}_{\nu,n}},
\end{equation}
respectively, where $\nu>-1$ and $m\in \mathbb{N}$. Note that for $m=1$, we have \cite{bpogsz}
\begin{equation}\label{rf3}
\eta_{2}(\nu)=\sum_{n\geq 1}\frac{1}{\alpha^2_{\nu,n}}=\frac{3}{4(\nu+1)}
\end{equation}
By using the series \eqref{sum}, one gets
$$
\frac{1}{4}\left(2\nu+1-\frac{zW'_{\nu}(z)}{W_{\nu}(z)}\right)=\sum_{n\geq 1}\frac{z^4}{\gamma^4_{\nu,n}-z^4},
$$
and taking limit $z\rightarrow 0$ followed by dividing with $z^4$ on each side, we obtain
\begin{equation}\label{rf4}
\zeta_{4}(\nu)=\sum_{n\geq 1}\frac{1}{\gamma^4_{\nu,n}}=\frac{1}{16(\nu+1)(\nu+2)(\nu+3)}=\frac{1}{2^4(\nu+1)_3}.
\end{equation}
For sake of brevity we denote $(\nu+1)(\nu+2)(\nu+3)$ by $(\nu+1)_3$ using the well known Pochhammer (Appel) symbol defined by $(\alpha)_0=0$ for $\alpha\neq 0$ and
$$
(\alpha)_n=\alpha(\alpha+1)(\alpha+2)\cdots(\alpha+n-1),~~~~\mbox{for}~~~~n\geq 1.
$$
In general, for any $m\in \mathbb{N}$, the Rayleigh function $\zeta_{4m}(\nu)$ can be obtained by comparing the coefficients of $z^{4m}$ on both sides of \eqref{series2}. For example, by comparing the coefficients of $z^{4}$ on both sides of \eqref{series2}, one can get \eqref{rf4} and by comparing the coefficients of $z^{8}$ on both sides of \eqref{series2} yields
\begin{equation}\label{rf5}
\zeta_{8}(\nu)=\frac{5\nu+17}{256(\nu+1)^2(\nu+2)^2(\nu+3)^2(\nu+4)(\nu+5)}
=\frac{5\nu+17}{2^8(\nu+1)_3(\nu+1)_5}.
\end{equation}
Alternatively, taking into account the power series \eqref{sum} and infinite product representation \eqref{product} one can extract the Rayleigh function $\zeta_{4m}(\nu)$ by using the Euler-Rayleigh method (see \cite[p. 3]{ismail}). Namely, let $f(z)$ be an entire function with power series representation
$$
f(z)=1+\sum_{n\geq1}a_nz^n
$$
and an infinite product representation
$$
f(z)=\prod_{n\geq1}\left(1-\frac{z}{z_n}\right),
$$
where it is assumed that $\sum\limits_{n\geq 1}|z_n|^{-1}<\infty$. Then the Rayleigh function
$$
S_m=\sum_{m\geq1}\frac{1}{z^m_k}
$$
is given by the following formula
$$
S_n=-na_n-\sum_{i=1}^{n-1}a_iS_{n-i}.
$$
Therefore, by taking $f(z)=\mathscr{W}_{\nu}(\sqrt[4]{z})$, from \eqref{sum} and \eqref{product} we have
$$
a_n=\frac{(-1)^n\Gamma(\nu+1)\Gamma(\nu+2)}{2^{4n}n!\Gamma(\nu+n+1)\Gamma(\nu+2n+2)}
$$
and hence $S_1=\zeta_{4}(\nu)=-a_1$, $S_2=\zeta_{8}(\nu)=-2a_2-a_1S_1=-2a_2+a^2_1$ and so on.

Now, we present the Euler-Rayleigh inequalities for zeros of Dini functions and cross-product of Bessel functions, which will be used in sequel.
\begin{lemma}
Let $\nu>-1$ and $m\in \mathbb{N}$. Then
\begin{equation}\label{ray-ineq1}
\left[\eta_{2m}(\nu)\right]^{-1/m}<\alpha^2_{\nu,1}<\frac{\eta_{2m}(\nu)}{\eta_{2m+2}(\nu)}
\end{equation}
and
\begin{equation}\label{ray-ineq2}
\left[\zeta_{4m}(\nu)\right]^{-1/m}<\gamma^4_{\nu,1}<\frac{\zeta_{4m}(\nu)}{\zeta_{4m+4}(\nu)}.
\end{equation}
\end{lemma}

The above inequalities can be verified easily by using the definition of $\eta_{2m}(\nu)$,
$\zeta_{4m}(\nu)$ and the order relations $0<\alpha_{\nu,1}<\alpha_{\nu,2}<\cdots<\alpha_{\nu,n}<\cdots$ and $0<\gamma_{\nu,1}<\gamma_{\nu,2}<\cdots<\gamma_{\nu,n}<\cdots$.

An immediate consequence of the above inequality will give the lower and upper bounds for the smallest positive zero of the cross-product of Bessel functions.
\begin{theorem}\label{theorem4'}
Let $\nu>-1$ and $\gamma_{\nu,1}$ denote the smallest positive zero of the cross-product of Bessel functions $W_{\nu}(z)$. Then we have the following bounds:
\begin{equation}\label{ub1}
\frac{2^4(\nu+1)_3\sqrt{(\nu+4)(\nu+5)}}{\sqrt{5\nu+17}}<\gamma^4_{\nu,1}<\frac{2^4(\nu+1)_5}{5\nu+17}.
\end{equation}
\end{theorem}

Note that using \eqref{rf4} and \eqref{rf5},  the left-hand side of the inequality \eqref{ub1} follows from left-hand side of \eqref{ray-ineq2} by taking $m=2,$ while the right-hand side of the inequality \eqref{ub1} can be extracted from right-hand side of \eqref{ray-ineq2} by taking $m=1$. So we omit the proof of Theorem \ref{theorem4'}.

Observe that for $m=1$, the left-hand side of \eqref{ray-ineq2} gives the inequality
$$
\gamma^4_{\nu,1}>2^4(\nu+1)_3,
$$
which is weaker than the left-hand side of \eqref{ub1}.

The power series representation \cite{bpogsz}
$$
\frac{zd'_{\nu}(z)}{d_{\nu}(z)}=\nu-2\sum_{m\geq 1}\eta_{2m}(\nu)z^{2m}
$$
which is valid for $\nu>-1$, $z\in \mathbb{C}$ such that $|z|<\alpha_{\nu,1}$, can be rewritten as
\begin{equation}\label{series1}
\frac{z\mathcal{D}'_{\nu}(z)}{\mathcal{D}_{\nu}(z)}=-2\sum_{m\geq 1}\eta_{2m}(\nu)z^{2m},
\end{equation}
where $\mathcal{D}_{\nu}(z)=2^{\nu}\Gamma(\nu+1)z^{-\nu}d_{\nu}(z)$. Therefore, the function $x\mapsto -\frac{x\mathcal{D}'_{\nu}(x)}{\mathcal{D}_{\nu}(x)}$ is absolutely monotonic on $(0,\alpha_{\nu,1})$ for all $\nu>-1$. The next theorem is analogous to this result.

\begin{theorem}\label{theorem5}
Let $\nu>-1$ and $z\in \mathbb{C}$ such that $|z|<\gamma_{\nu,1}$. Then
\begin{equation}\label{series2}
\frac{z\mathscr{W}'_{\nu}(z)}{\mathscr{W}_{\nu}(z)}=-4\sum_{m\geq 1}\zeta_{4m}(\nu)z^{4m}.
\end{equation}
Moreover, the function
$$
x\mapsto -\frac{x\mathscr{W}'_{\nu}(x)}{\mathscr{W}_{\nu}(x)}
$$
is absolutely monotonic on $(0,\gamma_{\nu,1})$ for all $\nu>-1$.
\end{theorem}
In addition, the next result is valid.
\begin{theorem}\label{theorem6}
Let $\mu\geq \nu>-1$. Then the functions $f_{\mu,\nu},$ $g_{\mu,\nu},$ $h_{\nu},$ $q_{\nu}:[0,\gamma^4_{\nu,1})\rightarrow (0,\infty)$, defined by
$$
f_{\mu,\nu}(x)=\left[\log\left(x^{\frac{\nu-\mu}{2}}e^{\frac{x}{16}
\left(\frac{1}{(\mu+1)_3}-\frac{1}{(\nu+1)_3}\right)}\frac{W_{\mu}(\sqrt[4]{x})}{W_{\nu}(\sqrt[4]{x})}
\right)\right]',
$$
$$
g_{\mu,\nu}(x)=x^{\frac{\nu-\mu}{2}}e^{\frac{x}{16}
\left(\frac{1}{(\mu+1)_3}-\frac{1}{(\nu+1)_3}\right)}\frac{W_{\mu}(\sqrt[4]{x})}{W_{\nu}(\sqrt[4]{x})},
$$
$$
h_{\nu}(x)=\left[\log \left(\frac{x^{\frac{\nu}{2}+\frac{1}{4}}e^{\frac{-x}{16(\nu+1)_3}}}{W_{\nu}(\sqrt[4]{x})}\right)\right]'
$$
and
$$
q_{\nu}(x)=\frac{x^{\frac{\nu}{2}+\frac{1}{4}}e^{\frac{-x}{16(\nu+1)_3}}}{W_{\nu}(\sqrt[4]{x})}
$$
are absolutely monotonic.
\end{theorem}
Observe that the above absolutely monotonicity of $q_{\nu}$ can be used to find the upper bound for the cross-product of Bessel and modified Bessel functions. Namely, we have the following inequality.
\begin{corollary}\label{corol6.1}
If $\nu>-1$ and $x\in [0,\gamma_{\nu,1})$, then
$$
W_{\nu}(x)\leq \frac{x^{2\nu+1}e^{-\frac{x^4}{16(\nu+1)_3}}}{2^{2\nu}\Gamma(\nu+1)\Gamma(\nu+2)}.
$$
\end{corollary}

\subsection{\bf Redheffer-type inequalities}
We continue with another set of results, namely Redheffer-type inequalities. In the literature the inequality
$$
\frac{\sin x}{x}\geq \frac{\pi^2-x^2}{\pi^2+x^2},~~~~\mbox{where}~~x\in\mathbb{R},
$$
is known as Redheffer inequality, see \cite{redheffer}. In \cite{baricz1} the author extended the above Redheffer type inequalities for the normalized Bessel functions of the first kind $\mathcal{J}_{\nu}(x)=2^{\nu}\Gamma(\nu+1)x^{-\nu}J_{\nu}(x)$ and the normalized modified Bessel function $\mathcal{I}_{\nu}(x)=2^{\nu}\Gamma(\nu+1)x^{-\nu}I_{\nu}(x).$ For more details about Redheffer type inequalities one can refer to \cite{baricz-wu,chen-zhao-qi,zhu1} and to the references therein. Recently in \cite{bps2}, Redheffer-type inequalities for modified Dini functions were studied.
In this subsection, we study Redheffer type inequalities for Dini functions and cross-product of Bessel and modified Bessel functions. Motivated by the result from \cite[Theorem 1]{zhu1}, we extend and sharpen the Redheffer-type inequalities for modified Dini functions \cite[Theorem 7]{bps2}.

\begin{theorem}\label{theorem7}
Let $\alpha_{\nu,n}$ and $\gamma_{\nu,n}$ denote the $n$-th positive zero of $d_{\nu}$ and $W_{\nu},$ respectively. The following Redheffer-type inequalities are valid:
\begin{enumerate}
\item[\bf a.] If $\nu>-1$ and $\Psi_{\nu}(n)=\alpha^2_{\nu,n+1}-\alpha_{\nu,1}\alpha_{\nu,n}-\alpha_{\nu,n}\alpha_{\nu,n+1}\geq 0$ for $n\in \mathbb{N},$ then
\begin{equation}\label{redheffer1}
\mathcal{D}_{\nu}(x)\geq \frac{\alpha^2_{\nu,1}-x^2}{\alpha^2_{\nu,1}+x^2}~~\mbox{for all}~~|x|\leq \delta_{\nu}=\min_{n\geq 1,\nu>-1}\left\{\alpha_{\nu,1},\sqrt{\Psi_{\nu}(n)}\right\}.
\end{equation}
\item[\bf b.] If $\nu\in (-1,8),$ then
\begin{equation}\label{redheffer2}
\mathcal{D}_{\nu}(x)\leq \left(\frac{\alpha^2_{\nu,1}-x^2}{\alpha^2_{\nu,1}+x^2}\right)^{m_{\nu}}~~\mbox{for all}~~x\in (-\alpha_{\nu,1},\alpha_{\nu,1}),
\end{equation}
where $m_{\nu}=\frac{3\alpha^2_{\nu,1}}{8(\nu+1)}$ is the best possible constant.
\item[\bf c.] If $\nu>-1$ and $\Omega_{\nu}(n)=\gamma^4_{\nu,n+1}-\gamma^2_{\nu,1}\gamma^2_{\nu,n}
    -\gamma^2_{\nu,n}\gamma^2_{\nu,n+1}
\geq 0$ for $n\in \mathbb{N},$ then
\begin{equation}\label{redheffer3}
\mathscr{W}_{\nu}(x)\geq \frac{\gamma^4_{\nu,1}-x^4}{\gamma^4_{\nu,1}+x^4}~~\mbox{for all}~~|x|\leq \epsilon_{\nu}=\min_{n\geq 1,\nu>-1}\left\{\gamma_{\nu,1},\sqrt{\Omega_{\nu}(n)}\right\}.
\end{equation}
\item[\bf d.] If $\nu\in (-1,r),$ where $r=\frac{1+\sqrt{57}}{2}$ is the positive root of $\nu-\nu^2+14=0,$ then
\begin{equation}\label{redheffer2'}
\mathscr{W}_{\nu}(x)\leq \left(\frac{\gamma^4_{\nu,1}-x^4}{\gamma^4_{\nu,1}+x^4}\right)^{n_{\nu}}~~\mbox{for all}~~x\in (-\gamma_{\nu,1},\gamma_{\nu,1}),
\end{equation}
where $n_{\nu}=\frac{\gamma^4_{\nu,1}}{32(\nu+1)_3}$ is the best possible constant.
\end{enumerate}
\end{theorem}
The corresponding result for the modified Dini functions reads as follows.
\begin{theorem}\label{theorem8}
Let $r\in (0,\infty)$, $|x|<r$, $\nu>-1$ and $\lambda_{\nu}$ be the modified Dini function defined by \eqref{mdfprod}. Then the following Redheffer-type inequality
\begin{equation}\label{redheffer4}
\left( \frac{r^2+x^2}{r^2-x^2}\right)^{\alpha}\leq \lambda_{\nu}(x)
\leq \left( \frac{r^2+x^2}{r^2-x^2}\right)^{\beta}
\end{equation}
holds if and only if $\alpha\leq 0$ and $\beta\geq \frac{3\alpha^2_{\nu,1}}{8(\nu+1)}$.
\end{theorem}

We would like to take the opportunity to correct a mistake in the paper \cite{baricz-wu}. In the final expression for $\varphi'_{\nu}(x)$ \cite[p. 263]{baricz-wu}, $\frac{1}{j^{2m-2}_{\nu,1}}$ should be replaced by $\frac{1}{j^{2m-4}_{\nu,1}},$ where $j_{\nu,1}$ stands for the first positive zero of the Bessel function $J_{\nu}.$ With this change, the following inequalities in \cite[p. 259]{baricz-wu} may not hold true for all $\nu\geq-7/8:$
\begin{equation}\label{redheffer5}
\mathcal{J}_{\nu}(x)\leq \left(\frac{j^2_{\nu,1}-x^2}{j^2_{\nu,1}+x^2}\right)^{\beta_{\nu}} ~~\mbox{for all}~~|x|< j_{\nu,1}
\end{equation}
and
\begin{equation}\label{redheffer6}
\frac{\mathcal{J}_{\nu+1}(x)}{\mathcal{J}_{\nu}(x)}\geq \left(\frac{j^2_{\nu,1}+x^2}{j^2_{\nu,1}-x^2}\right)^{\gamma_{\nu}} ~~\mbox{for all}~~|x|< j_{\nu,1},
\end{equation}
where $\beta_{\nu}=\frac{j^2_{\nu,1}}{8(\nu+1)}$ and $\gamma_{\nu}=\frac{j^2_{\nu,1}}{8(\nu+1)(\nu+2)}$ are best possible constants. Nevertheless, the above inequalities are valid for $\nu\in(-1,\nu_0]$ and $|x|<j_{\nu,1}$, where $\nu_0\in(1,2)$ is the unique root of the equation $j^2_{\nu,1}= 8(\nu+1)$. Before we prove the above inequalities, let us recall the following \cite[Lemma 1]{baricz-wu}, which will be useful in the sequel.
\begin{lemma}\label{lemma2}
Let $\nu>-1$ and $j_{\nu,1}$ be the first positive zero of the Bessel function $J_{\nu}$. Then the equation $j^2_{\nu,1}= 8(\nu+1)$  has exactly one positive root $\nu_0\in(1,2)$. Moreover,
$$
\begin{cases}
~   j^2_{\nu,1}\leq 8(\nu+1) &   \mbox{for $\nu\in (-1,\nu_0]$,}\\
~  j^2_{\nu,1}\geq 8(\nu+1) &   \mbox{for $\nu\geq \nu_0$}.
\end{cases}
$$
\end{lemma}

Now, taking into account the above correction in the expression $\varphi'_{\nu}(x)$ \cite{baricz-wu} we have,
\begin{eqnarray*}
\varphi'_{\nu}(x)&=&\frac{1}{2(\nu+1)x}\frac{1}{j^2_{\nu,1}+x^2}\sum_{m\geq 2}\left[4(\nu+1)j^2_{\nu,1}\sigma^{(2m)}_{\nu}
+4(\nu+1)\sigma^{(2m-2)}_{\nu}-\frac{1}{j^{2m-4}_{\nu,1}}\right]x^{2m}\\
&\geq &\frac{1}{2(\nu+1)x}\frac{1}{j^2_{\nu,1}+x^2}\sum_{m\geq 2}\left[8(\nu+1)j^2_{\nu,1}\sigma^{(2m)}_{\nu}-\frac{1}{j^{2m-4}_{\nu,1}}\right]x^{2m}\\
&=&\frac{1}{2(\nu+1)x}\frac{1}{j^2_{\nu,1}+x^2}\sum_{m\geq 2}\left[\frac{8(\nu+1)j^{-2}_{\nu,1}j^{2m}_{\nu,1}\sigma^{(2m)}_{\nu}-1}{j^{2m-4}_{\nu,1}}\right]x^{2m}\\
&\geq &\frac{1}{2(\nu+1)x}\frac{1}{j^2_{\nu,1}+x^2}\sum_{m\geq 2}\left[\frac{8(\nu+1)j^{-2}_{\nu,1}-1}{j^{2m-4}_{\nu,1}}\right]x^{2m},
\end{eqnarray*}
which in view of the Lemma \ref{lemma2} gives
$$
\varphi'_{\nu}(x)\geq 0~~\mbox{for}~~\nu\in(-1,\nu_0].
$$
Here we have used the Euler-Rayleigh inequalities \cite[p, 502]{watson}
\begin{equation}\label{ray-ineq3}
\left(\sigma^{(2m)}_{\nu}\right)^{-1/m}<j^2_{\nu,1}<\frac{\sigma^{(2m)}_{\nu}}{\sigma^{(2m+2)}_{\nu}},
\end{equation}
which are valid for $m\in \mathbb{N}$ and $\nu>-1$, where
$$
\sigma^{(2m)}_{\nu}=\sum_{n\geq 1}\frac{1}{j^{2m}_{\nu,n}}
$$
is the Rayleigh function of order $2m$. The rest of the proof is same as in \cite[p. 263]{baricz-wu}.

It is also interesting to note that for $\nu\geq \nu_0$, the following new Redheffer type inequalities hold.
\begin{theorem}\label{theorem9}
Let $\nu\geq \nu_0$, where $\nu_0\in (1,2)$ is the unique root of the equation $j^2_{\nu,1}= 8(\nu+1)$. Then the following new Redheffer type inequalities are valid
\begin{equation}\label{redheffer7}
\mathcal{J}_{\nu}(x)\leq \frac{j^2_{\nu,1}-x^2}{j^2_{\nu,1}+x^2} ~~~~\mbox{for all}~~|x|< j_{\nu,1},
\end{equation}
and
\begin{equation}\label{redheffer8}
\frac{\mathcal{J}_{\nu+1}(x)}{\mathcal{J}_{\nu}(x)}\geq \left(\frac{j^2_{\nu,1}+x^2}{j^2_{\nu,1}-x^2}\right)^{\frac{1}{\nu+2}} ~~\mbox{for all}~~|x|< j_{\nu,1}.
\end{equation}
\end{theorem}

\subsection{Bounds for logarithmic derivative of Bessel related functions} In this section we investigate the bounds for logarithmic derivative of Dini functions and the logarithmic derivative of cross-product of Bessel and modified Bessel functions. The idea of these results come from \cite{zhu2}.

\begin{theorem}\label{theorem10}
Let $\nu>-1$, $A_n=\alpha^2_{\nu,1}\eta_{2n+2}(\nu)-\eta_{2n}(\nu)$, $n\in \mathbb{N}$. Then for $n\in \mathbb{N},$ the following inequality holds true for all $0<|x|<\alpha_{\nu,1},$
$$
\frac{R_{2n}(x)+\frac{4}{3}(\nu+1)ax^{2n+2}}{\alpha^2_{\nu,1}-x^2}<-\frac{2(\nu+1)}{3x}\cdot\frac{\mathcal{D}'_{\nu}(x)}{\mathcal{D}_{\nu}(x)}
<\frac{R_{2n}(x)+\frac{4}{3}(\nu+1)bx^{2n+2}}{\alpha^2_{\nu,1}-x^2},
$$
where
$$a=\frac{1}{\alpha^{2n+2}_{\nu,1}}\left(1-\frac{3\alpha^2_{\nu,1}}{4(\nu+1)}-\sum_{m=1}^nA_m\alpha^{2m}_{\nu,1}\right),~~ b=A_{n+1}~~ \mbox{and}~~ R_{2n}(x)=\alpha^2_{\nu,1}+\frac{4(\nu+1)}{3}\sum_{m=1}^nA_mx^{2m}.$$
Moreover, $a$ and $b$ are sharp.
\end{theorem}

\begin{theorem}\label{theorem11}
Let $\nu>-1$, $B_n=\gamma^4_{\nu,1}\zeta_{4n+4}(\nu)-\zeta_{4n}(\nu)$, $n\in \mathbb{N}$. Then for $n\in \mathbb{N},$ the following inequality holds true for all $0<|x|<\gamma_{\nu,1},$
$$
\frac{S_{4n}(x)+16(\nu+1)_3rx^{4n+4}}{\gamma^4_{\nu,1}-x^4}<-\frac{4(\nu+1)_3}{x^3}\cdot\frac{\mathscr{W}'_{\nu}(x)}{\mathscr{W}_{\nu}(x)}
<\frac{S_{4n}(x)+16(\nu+1)_3sx^{4n+4}}{\gamma^4_{\nu,1}-x^4},
$$
where
$$r=\frac{1}{\gamma^{4n+4}_{\nu,1}}\left(1-\frac{\gamma^4_{\nu,1}}{16(\nu+1)_3}-\sum_{m=1}^nB_m\gamma^{4m}_{\nu,1}\right), ~~s=B_{n+1}~~ \mbox{and}~~ S_{4n}(x)=\gamma^4_{\nu,1}+16(\nu+1)_3\sum_{m=1}^nB_mx^{4m}.$$
Moreover, $r$ and $s$ are sharp.
\end{theorem}

\section{\bf Proofs of main results}\label{proof}
In this section we prove our main results.
\begin{proof}[\bf Proof of Theorem \ref{theorem1}]
{\bf a.} By using the infinite product representation \eqref{product} and the order relation
$$
0<\gamma_{\nu,1}<\gamma_{\nu,2}<\cdots<\gamma_{\nu,n}<\cdots,
$$
we note that if $x\in [\gamma_{\nu,2n-1},\gamma_{\nu,2n}]$ or $x\in [-\gamma_{\nu,2n},-\gamma_{\nu,2n-1}]$ then the first $(2n-1)$
terms of the product \eqref{product} are negative and the remaining terms are strictly positive. Therefore $\mathscr{W}_{\nu}(x)$ becomes negative on $S$. Now if $x\in (-\gamma_{\nu,1},\gamma_{\nu,1})$ then each terms of the product \eqref{product}
are strictly positive and if $x\in (\gamma_{\nu,2n},\gamma_{\nu,2n+1})$ or $x\in (-\gamma_{\nu,2n+1},-\gamma_{\nu,2n})$, then
the first $2n$ terms are strictly negative while the remaining terms are strictly positive. Therefore $\mathscr{W}_{\nu}(x)>0$
on $\mathbb{R}\setminus S$.

{\bf b.} From part {\bf a}, we have $\mathscr{W}_{\nu}(x)>0$ for $x\in (-\gamma_{\nu,1},\gamma_{\nu,1})$. Therefore the infinite product representation \eqref{product} gives
$$
\left(\log\mathscr{W}_{\nu}(x)\right)'=\frac{\mathscr{W}'_{\nu}(x)}{\mathscr{W}_{\nu}(x)}=-\sum_{n\geq 1}\frac{4x^3}{\gamma^4_{\nu,n}-x^4}
$$
and hence the function $x\mapsto\mathscr{W}_{\nu}(x)$ is strictly increasing on $(-\gamma_{\nu,1},0]$ and strictly decreasing on $[0,\gamma_{\nu,1})$.

{\bf c.} By using the above equation and part {\bf a} of this theorem, we have
$$
\left(\log\mathscr{W}_{\nu}(x)\right)''=\left(\frac{\mathscr{W}'_{\nu}(x)}{\mathscr{W}_{\nu}(x)}\right)'=-\sum_{n\geq 1}\frac{4x^2(3\gamma^4_{\nu,n}+x^4)}{(\gamma^4_{\nu,n}-x^4)^2}
$$
and
$$
\left(\frac{x\mathscr{W}'_{\nu}(x)}{\mathscr{W}_{\nu}(x)}\right)'=-\sum_{n\geq 1}\frac{16x^3\gamma^4_{\nu,n}}{(\gamma^4_{\nu,n}-x^4)^2}.
$$
From this we conclude that the function $x\mapsto\mathscr{W}_{\nu}(x)$ is strictly log-concave on $\mathbb{R}\setminus S$ and
strictly geometrically concave on $(0,\infty)\setminus S_2$.

{\bf d.} Since the function $x\mapsto x^{2\nu+1}$ is log-concave on $(0,\infty)$ for all $\nu\geq -\frac{1}{2}$ and from part {\bf c} the function $x\mapsto\mathscr{W}_{\nu}(x)$ is strictly log-concave on $\mathbb{R}\setminus S$, we conclude that the functions
$$x\mapsto W_{\nu}(x)=\frac{x^{2\nu+1}\mathscr{W}_{\nu}(x)}{2^{2\nu}\Gamma(\nu+1)\Gamma(\nu+2)}$$
is strictly log-concave on $(0,\infty)\setminus S_2$ for all $\nu\geq-\frac{1}{2}$.
Here we used the fact that product of a log-concave function and a strictly log-concave function is strictly log-concave.

{\bf e.} By using again the infinite product representation \eqref{product}, we get
$$
\frac{\partial}{\partial\nu}\left(\log\mathscr{W}_{\nu}(x)\right)=\sum_{n\geq 1}\frac{4x^4\frac{\partial}{\partial\nu}(\gamma_{\nu,n})}
{\gamma_{\nu,n}(\gamma^4_{\nu,n}-x^4)}
$$
and
$$\frac{\partial}{\partial\nu}\left(\frac{x\mathscr{W}'_{\nu}(x)}
{\mathscr{W}_{\nu}(x)}\right) =\sum_{n\geq 1}\frac{16x^4\gamma^3_{\nu,n}\frac{\partial\gamma_{\nu,n}}
{\partial\nu}}{(\gamma^4_{\nu,n}-x^4)^2}.
$$
From these expressions and the result \cite[Lemma 4]{abp}: $\nu\mapsto \gamma_{\nu,n}$ is increasing on $(-1,\infty)$,  the desired conclusion follows.

{\bf f.} From the infinite product representation \eqref{product}, we have
$$
\left(-\log\mathscr{W}_{\nu}(\sqrt[4]{x})\right)'=\sum_{n\geq 1}\frac{1}{\gamma^4_{\nu,n}-x},
$$
which is absolutely monotonic on $(0,\gamma^4_{\nu,1})$ for all $\nu>-1$. Since the exponential of a function having an absolutely monotonic derivative is absolutely monotonic, we conclude that
the function $x\mapsto1/\mathscr{W}_{\nu}(\sqrt[4]{x})$ is absolutely monotonic on $(0,\gamma^4_{\nu,1})$ for all $\nu>-1$.
\end{proof}

\begin{figure}[!ht]
   \centering
       \includegraphics[width=12cm]{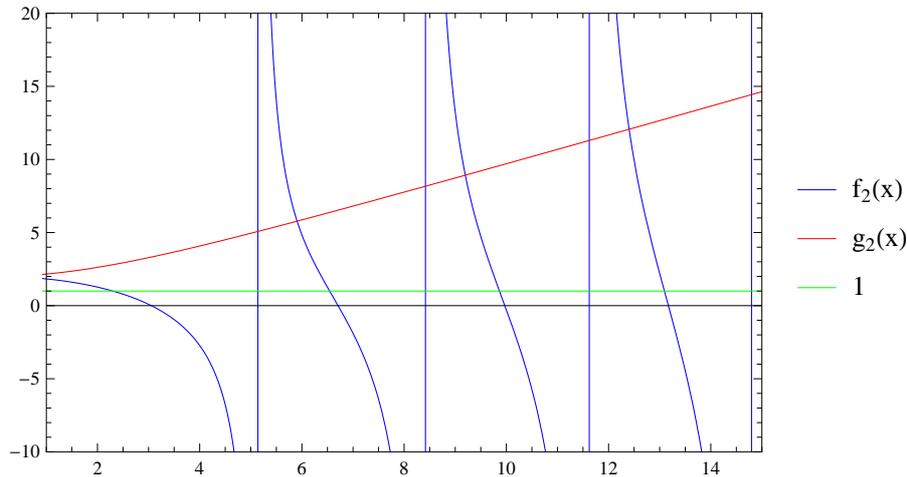}
       \caption{Interlacing of zeros of Dini functions and cross-product of
Bessel functions: the graph of the functions $x\mapsto$ f$_2(x)=\frac{xJ_2'(x)}{J_2(x)}$ and $x\mapsto
$g$_2(x)=\frac{xI_2'(x)}{I_2(x)}$ on $[0,15]$.}
       \label{fig1}
\end{figure}

\begin{proof}[\bf Proof of Theorem \ref{theorem2}]
By using the inequalities \eqref{IP1} and \eqref{IP2}, we have
$$
j_{\nu,n-1}<\alpha_{\nu,n}<j_{\nu,n}<\gamma_{\nu,n}<j_{\nu,n+1},
$$
where $ n\in \mathbb{N}$ and hence the left-hand side of the inequality \eqref{IP3} follows. To prove the right-hand side of the inequality \eqref{IP3}, observe that the zeros of the cross-product of Bessel functions \eqref{sum} are the roots of the equation
\begin{equation}\label{cpzero}
\frac{zJ'_{\nu}(z)}{J_{\nu}(z)}=\frac{zI'_{\nu}(z)}{I_{\nu}(z)}
\end{equation}
and the zeros of Dini function $z\mapsto (1-\nu)J_{\nu}(z)+zJ'_{\nu}(z)$ are roots of the equation
\begin{equation}\label{dinizero}
\frac{zJ'_{\nu}(z)}{J_{\nu}(z)}=\nu-1.
\end{equation}
Now in view of the infinite product representations of Bessel and modified Bessel functions of first kind, namely,
$$
2^{\nu}\Gamma(\nu+1)x^{-\nu}J_{\nu}(x)=\prod_{n\geq 1}\left(1-\frac{x^2}{j^2_{\nu,n}}\right)\ \ \ \mbox{and}\ \ \
2^{\nu}\Gamma(\nu+1)x^{-\nu}I_{\nu}(x)=\prod_{n\geq 1}\left(1+\frac{x^2}{j^2_{\nu,n}}\right),
$$
we obtain
$$
\left(\frac{xJ'_{\nu}(x)}{J_{\nu}(x)}\right)'=-\sum_{n\geq 1}\frac{4xj^2_{\nu,n}}{(j^2_{\nu,n}-x^2)^2}\ \ \ \mbox{and}\ \ \
\left(\frac{xI'_{\nu}(x)}{I_{\nu}(x)}\right)'=\sum_{n\geq 1}\frac{4xj^2_{\nu,n}}{(j^2_{\nu,n}+x^2)^2},
$$
respectively. Therefore, for $\nu>-1$ the function $x\mapsto \frac{xJ'_{\nu}(x)}{J_{\nu}(x)}$ is strictly decreasing on each interval $(j_{\nu,n},j_{\nu,n+1})$, $n\in \mathbb{N}$ and the function $x\mapsto \frac{xI'_{\nu}(x)}{I_{\nu}(x)}$ is strictly increasing on $(0,\infty)$.
This implies that there exists a unique root $\gamma_{\nu,n}$ of the equation \eqref{cpzero} and a unique root $\alpha_{\nu,n+1}$ of the equation \eqref{dinizero} in each interval $(j_{\nu,n},j_{\nu,n+1})$ for all $n\in \mathbb{N}$. Since the function $x\mapsto \frac{xJ'_{\nu}(x)}{J_{\nu}(x)}$ is strictly decreasing on each interval $(j_{\nu,n},j_{\nu,n+1})$, $n\in \mathbb{N}$ and we have the limit
$$
\lim_{x\rightarrow 0}\frac{xI'_{\nu}(x)}{I_{\nu}(x)}=\nu>\nu-1,
$$
we conclude that $\gamma_{\nu,n}<\alpha_{\nu,n+1}$ for all $n\in \mathbb{N}$.  This interlacing property is illustrated in Figure \ref{fig1} for $\nu=2$ and $x\in(0,j_{2,4})$. This completes the proof.
\end{proof}

\begin{proof}[\bf Proof of Corollary \ref{corol2.1}]
Using \eqref{IP3}, we have for all $x\in (-\gamma_{\nu,1},\gamma_{\nu,1})$
$$\prod_{n\geq 1}\left(1-\frac{x^4}{\alpha^4_{\nu,n}}\right) < \prod_{n\geq 1}\left(1-\frac{x^4}{\gamma_{\nu,n}^4}\right)
< \prod_{n\geq 1}\left(1-\frac{x^4}{\alpha_{\nu,n+1}^4}\right),
$$
which on using \eqref{product}, \eqref{dfprod} and \eqref{mdfprod} gives the inequality \eqref{corol2.1_ine1}. Similarly, by using the interlacing inequality \eqref{IP2} one can extract the inequality \eqref{corol2.1_ine2}. To prove the inequality \eqref{corol2.1_ine3}, observe that the inequality \eqref{corol2.1_ine2} can be rewritten as
$$
1<\frac{(\nu+1)}{x}\left(\frac{I_{\nu+1}(x)}{I_{\nu}(x)}+\frac{J_{\nu+1}(x)}{J_{\nu}(x)}\right)<\frac{j^4_{\nu,1}}{j^4_{\nu,1}-x^4} ~~~~\mbox{for}~~~~|x|<j_{\nu,1},
$$
which in view of the formulas
$$\mathcal{I}'_{\nu}(x)=2^{\nu}\Gamma(\nu+1)(x^{-\nu}I_{\nu}(x))'
=2^{\nu}\Gamma(\nu+1)x^{-\nu}I_{\nu+1}(x)$$
and
$$\mathcal{J}'_{\nu}(x)=2^{\nu}\Gamma(\nu+1)(x^{-\nu}J_{\nu}(x))' =-2^{\nu}\Gamma(\nu+1)x^{-\nu}J_{\nu+1}(x),$$
is equivalent to
\begin{equation}\label{ratio_bound1}
1<\frac{(\nu+1)}{x}\left(\frac{\mathcal{I}'_{\nu}(x)}{\mathcal{I}_{\nu}(x)}-\frac{\mathcal{J}'_{\nu}(x)}{\mathcal{J}_{\nu}(x)}\right)<\frac{j^4_{\nu,1}}{j^4_{\nu,1}-x^4} ~~~~\mbox{for}~~~~|x|<j_{\nu,1},
\end{equation}
Now for $x\in (0,j_{\nu,1})$, integrating \eqref{ratio_bound1} we obtain
$$
\int_0^x\frac{t}{(\nu+1)}dt<\int_0^x\left(\frac{\mathcal{I}'_{\nu}(t)}{\mathcal{I}_{\nu}(t)}-\frac{\mathcal{J}'_{\nu}(t)}{\mathcal{J}_{\nu}(t)}\right)dt<\frac{j^4_{\nu,1}}{(\nu+1)}\int_0^x\frac{t}{j^4_{\nu,1}-t^4}dt,
$$
which implies that
$$
e^{\frac{x^2}{2(\nu+1)}}<\frac{\mathcal{I}_{\nu}(x)}{\mathcal{J}_{\nu}(x)}<\left(\frac{j^2_{\nu,1}+x^2}{j^2_{\nu,1}-x^2}\right)^{\frac{j^2_{\nu,1}}{4(\nu+1)}}.
$$
This proves the inequality \eqref{corol2.1_ine3}.
\end{proof}
\begin{proof}[\bf Proof of Theorem \ref{theorem3}]
{\bf a \& b.}  The normalized Dini function $z\mapsto \mathcal{D}_{\nu}(z)$ and the cross-product of Bessel and modified Bessel functions $z\mapsto\mathscr{W}_{\nu}(\sqrt{z})$ are entire functions of order $1/2$ and $1/4,$ respectively (see \cite{bpogsz}, \cite{abp}). Therefore the genus of the entire functions $z\mapsto\mathcal{D}_{\nu}(z)$ and $z\mapsto\mathscr{W}_{\nu}(\sqrt{z})$ is $0,$ as the genus of entire function of order $\rho$ is $[\rho]$ when $\rho$ is not an integer \cite[p. 34]{boas}. We also note that the zeros of $z\mapsto\mathcal{D}_{\nu}(z)$ and $z\mapsto\mathscr{W}_{\nu}(\sqrt{z})$ are all real when $\nu>-1$. Now recall Laguerre's theorem on separation of zeros \cite[p. 23]{boas} which states that, if $z\mapsto f(z)$ is a non-constant entire function, which is real for real $z$ and has only real zeros, and is of genus $0$ or $1$, then the zeros of $f'$ are also real and separated by the zeros of $f$. Therefore in view of Laguerre's theorem the conclusions follow.

{\bf c.} Since for $\nu>-1$ the function $\mathcal{D}_{\nu}$ belongs to Laguerre-P\'olya class of entire functions, it satisfies the Laguerre inequality \cite{skov}
$$
\left[\mathcal{D}^{(m)}_{\nu}(x)\right]^2-\mathcal{D}^{(m-1)}_{\nu}(x)\mathcal{D}^{(m+1)}_{\nu}(x)\geq 0.
$$
Using the derivative formulas
\begin{equation}\label{derivative1}
\mathcal{D}'_{\nu}(x)=2^{\nu}\Gamma(\nu+1)x^{-\nu-1}\left[xd'_{\nu}(x)-\nu d_{\nu}(x)\right]
\end{equation}
and
\begin{equation}\label{derivative2}
\mathcal{D}''_{\nu}(x)=2^{\nu}\Gamma(\nu+1)x^{-\nu-2}\left[x^2d''_{\nu}(x)-2\nu x d'_{\nu}(x)+\nu(\nu+1)d_{\nu}(x)\right],
\end{equation}
the above inequality for $m=1$ is equivalent to
$$
2^{2\nu}\Gamma^2(\nu+1)x^{-2\nu-2}\left[x^2(d'_{\nu}(x))^2-\nu d^2_{\nu}(x)-x^2d_{\nu}(x)d''_{\nu}(x)\right]\ge 0
$$
which implies that
$$
(d'_{\nu}(x))^2-d_{\nu}(x)d''_{\nu}(x)\geq \frac{\nu}{x^2} d^2_{\nu}(x)>0
$$
for $\nu>0$ and $x\in \mathbb{R},$ $x\neq0.$ Therefore the function $x\mapsto \frac{d'_{\nu}(x)}{d_{\nu}(x)}$ is strictly decreasing on $(0,\infty)\setminus \{\alpha_{\nu,n}\mid n\in \mathbb{N}\}$. In view of \cite[Lemma 2.2]{bs}, all zeros of $d_{\nu}(x)$ are real and simple and hence
$d'_{\nu}(x)\neq 0$ at $x=\alpha_{\nu,n}$, $n\in \mathbb{N}$. Thus, for a fixed $n\in \mathbb{N}$, we have the limit $\lim_{x\searrow \alpha_{\nu,n-1}} \frac{d'_{\nu}(x)}{d_{\nu}(x)}=\infty$ and $\lim_{x\nearrow \alpha_{\nu,n}} \frac{d'_{\nu}(x)}{d_{\nu}(x)}=-\infty$. Since the function $x\mapsto \frac{d'_{\nu}(x)}{d_{\nu}(x)}$ is strictly decreasing on $(0,\infty)\setminus \{\alpha_{\nu,n}\mid n\in \mathbb{N}\}$ it follows that in each interval $(\alpha_{\nu,n-1},\alpha_{\nu,n})$ there exists a unique zero $\alpha'_{\nu,n}$ of $d'_{\nu}(x)$. Here we used the convention that $\alpha_{\nu,0}=0$.
\end{proof}

\begin{proof}[\bf Proof of Theorem \ref{theorem4}]
From the infinite product representations \eqref{product}, \eqref{dfprod} and \eqref{mdfprod} it is easy to verify that for all $\nu>-1$, the functions $\mathscr{W}_{\nu}$, $\mathcal{D}_{\nu}$ and $\lambda_{\nu}$ satisfy the following identities (in other words, Mittag-Leffler expansions)
\begin{equation}\label{logderivative}
\frac{\mathscr{W}'_{\nu}(x)}{\mathscr{W}_{\nu}(x)}=\sum_{n\geq 1}\frac{-4x^3}{\gamma^4_{\nu,n}-x^4},
\end{equation}
\begin{equation}\label{logderivative1}
\frac{\mathcal{D}'_{\nu}(x)}{\mathcal{D}_{\nu}(x)}=\sum_{n\geq 1}\frac{-2x}{\alpha_{\nu,n}^2-x^2},
\end{equation}
and
\begin{equation}\label{logderivative2}
\frac{\lambda'_{\nu}(x)}{\lambda_{\nu}(x)}=\sum_{n\geq 1}\frac{2x}{\alpha_{\nu,n}^2+x^2}.
\end{equation}
In view of the above logarithmic derivative \eqref{logderivative1} of Dini functions $\mathcal{D}_{\nu}$, we obtain that
\begin{eqnarray*}
\Theta_1&=&\sum_{n\ge 1, n\neq k}\frac{1}{\alpha^2_{\nu,n}-\alpha^2_{\nu,k}}=\lim_{x\rightarrow\alpha_{\nu,k}}
\left[-\frac{1}{2x}\cdot\frac{\mathcal{D}'_{\nu}(x)}{\mathcal{D}_{\nu}(x)}-\frac{1}{\alpha^2_{\nu,k}-x^2}\right]\\
&=&-\frac{1}{2\alpha_{\nu,k}}\lim_{x\rightarrow\alpha_{\nu,k}}
\left[\frac{\mathcal{D}'_{\nu}(x)(\alpha^2_{\nu,k}-x^2)+2x\mathcal{D}_{\nu}(x)}{\mathcal{D}_{\nu}(x)(\alpha^2_{\nu,k}-x^2)}\right].
\end{eqnarray*}
Now, by applying the Bernoulli-L'Hospital rule twice and using the derivative formulas \eqref{derivative1} and \eqref{derivative2} we have
$$
\Theta_1=-\frac{1}{4\alpha_{\nu,k}}\lim_{x\rightarrow \alpha_{\nu,k}}\left[\frac{d''_{\nu}(x)}{d'_{\nu}(x)}-\frac{(2\nu+1)}{x}\right].
$$
Using the differential equation \cite[p. 13]{erdelyi}
$$
x^2(x^2-2\nu+1)d''_{\nu}(x)-x(x^2+2\nu-1)d'_{\nu}(x)-\left[(x^2-\nu^2)(x^2-2\nu+1)+2(1-\nu)x^2\right]d_{\nu}(x)=0,
$$
satisfied by the Dini function $d_{\nu}$, we obtain
\begin{equation}\label{limit}
\lim_{x\rightarrow \alpha_{\nu,k}}\frac{d''_{\nu}(x)}{d'_{\nu}(x)}=\frac{\alpha^2_{\nu,k}+2\nu-1}{\alpha_{\nu,k}(\alpha^2_{\nu,k}-2\nu+1)},
\end{equation}
and hence
$$
\Theta_1=-\frac{1}{4\alpha_{\nu,k}}\left[\frac{\alpha^2_{\nu,k}+2\nu-1}{\alpha_{\nu,k}(\alpha^2_{\nu,k}-2\nu+1)}-\frac{2\nu+1}{\alpha_{\nu,k}}\right].
$$
Therefore the relation \eqref{id1} is indeed true.

To prove the identity \eqref{id2}, we appeal to the formulas \eqref{logderivative1} and \eqref{logderivative2} to obtain
\begin{eqnarray*}
\Theta_2&=&\sum_{n\ge 1, n\neq k}\frac{1}{\alpha^4_{\nu,n}-\alpha^4_{\nu,k}}=\lim_{x\rightarrow\alpha_{\nu,k}}
\left[-\frac{1}{4x^3}\cdot\left(\frac{\mathcal{D}'_{\nu}(x)}{\mathcal{D}_{\nu}(x)}+\frac{\lambda'_{\nu}(x}{\lambda_{\nu}(x)}\right)-\frac{1}{\alpha^4_{\nu,k}-x^4}\right]\\
&=&-\frac{1}{4\alpha^3_{\nu,k}}\frac{\lambda'_{\nu}(\alpha_{\nu,k})}{\lambda_{\nu}(\alpha_{\nu,k})}-\lim_{x\rightarrow\alpha_{\nu,k}}
\frac{1}{4x^3}\cdot\left(\frac{\mathcal{D}'_{\nu}(x)}{\mathcal{D}_{\nu}(x)}+\frac{4x^3}{\alpha^4_{\nu,k}-x^4}\right)\\
&=&-\frac{1}{4\alpha^3_{\nu,k}}\sum_{n\geq 1}\frac{2\alpha_{\nu,k}}{\alpha^2_{\nu,n}+\alpha^2_{\nu,k}}-\frac{1}{4\alpha^3_{\nu,k}}\lim_{x\rightarrow\alpha_{\nu,k}}
\left(\frac{\mathcal{D}'_{\nu}(x)(\alpha^4_{\nu,k}-x^4)+4x^3\mathcal{D}_{\nu}(x)}{\mathcal{D}_{\nu}(x)(\alpha^4_{\nu,k}-x^4)}\right).
\end{eqnarray*}
Now, by applying again the Bernoulli-L'Hospital rule twice and using the derivative formulas \eqref{derivative1} and \eqref{derivative2} we obtain
$$
\Theta_2=-\frac{1}{2\alpha^2_{\nu,k}}\sum_{n\geq 1}\frac{1}{\alpha^2_{\nu,n}+\alpha^2_{\nu,k}}
-\frac{1}{4\alpha^3_{\nu,k}}\lim_{x\rightarrow \alpha_{\nu,k}}\left[\frac{1}{2}\cdot\frac{d''_{\nu}(x)}{d'_{\nu}(x)}-\frac{(2\nu+3)}{2x}\right],
$$
which on using the limit \eqref{limit} gives \eqref{id2}.

To prove the identity \eqref{id3}, first we will show that for $\nu>-1$ and $z\in \mathbb{C}$ we have
\begin{equation}\label{product1}
\frac{2^{2\nu}\Gamma(\nu+1)\Gamma(\nu+2)}{(2\nu+1)}z^{-2\nu}W'_{\nu}(z)=\prod_{n\geq 1}\left(1-\frac{z^4}{\gamma'^4_{\nu,n}}\right).
\end{equation}
To deduce the above Hadamard factorization of $W_{\nu}$, it is enough to show
\begin{equation}\label{product2}
\frac{2^{2\nu}\Gamma(\nu+1)\Gamma(\nu+2)}{(2\nu+1)}z^{-\nu}W'_{\nu}(\sqrt{z})=\prod_{n\geq 1}\left(1-\frac{z^2}{\gamma'^4_{\nu,n}}\right).
\end{equation}
Now, by using the power series representation \eqref{sum}, we have
$$
\frac{2^{2\nu}\Gamma(\nu+1)\Gamma(\nu+2)}{(2\nu+1)}z^{-\nu}W'_{\nu}(\sqrt{z})=1+\sum_{n\geq 1}\frac{(-1)^n(2\nu+4n+1)\Gamma(\nu+1)\Gamma(\nu+2)z^{2n}}{n!\Gamma(\nu+n+1)\Gamma(\nu+2n+2)2^{4n}(2\nu+1)}.
$$
This is an entire function of growth order $\frac{1}{4},$ since
$$
\lim_{n\rightarrow \infty}\frac{n\log n}{\log \Gamma(n+1)+\log \Gamma(\nu+n+1)+\log \Gamma(\nu+2n+1)+\log \frac{2^{4n}(2\nu+1)}{\Gamma(\nu+1)\Gamma(\nu+2)}-\log (2\nu+4n+1)}=\frac{1}{4},
$$
where this limit follows easily on considering the limit
$$
\lim_{n\rightarrow \infty}\frac{\log \Gamma(an+b)}{n\log n}=a~~~~\mbox{where}~~~~a,b>0.
$$
By applying Hadamard's Theorem \cite[p. 26]{lev} it follows that \eqref{product2} is indeed valid and consequently we get \eqref{product1}.

Now, we use the formula \eqref{logderivative} and we get
\begin{eqnarray*}
\Theta_3&=&\sum_{n\ge 1, n\neq k}\frac{1}{\gamma^4_{\nu,n}-\gamma^4_{\nu,k}}=\lim_{x\rightarrow \gamma_{\nu,k}}\left[-\frac{1}{4x^3}\cdot\frac{\mathscr{W}'_{\nu}(x)}{\mathscr{W}_{\nu}(x)}-\frac{1}{\gamma^4_{\nu,k}-x^4}\right]\\
&=&-\frac{1}{4\gamma^3_{\nu,k}}\lim_{x\rightarrow \gamma_{\nu,k}}
\left(\frac{\mathscr{W}'_{\nu}(x)(\gamma^4_{\nu,k}-x^4)+4x^3\mathscr{W}_{\nu}(x)}{\mathscr{W}_{\nu}(x)(\gamma^4_{\nu,k}-x^4)}\right),
\end{eqnarray*}
which on applying the Bernoulli-L'Hospital rule two times gives
$$
\Theta_3=\frac{1}{8\gamma^3_{\nu,k}}\lim_{x\rightarrow \gamma_{\nu,k}}
\left[\frac{3}{x}-\frac{\mathscr{W}''_{\nu}(x)}{\mathscr{W}'_{\nu}(x)}\right].
$$
Logarithmic differentiation of \eqref{product1} gives
$$\frac{W''_{\nu}(x)}{W'_{\nu}(x)}=\frac{2\nu}{x}-\sum_{n\geq 1}\frac{4x^3}{\gamma'^4_{\nu,n}-x^4}.
$$
Now, using the following derivative formulas which follow easily from \eqref{product}, we obtain that
$$\mathscr{W}'_{\nu}(x)=2^{\nu}\Gamma(\nu+1)\Gamma(\nu+2)x^{-2\nu-2}\left[xW'_{\nu}(x)-(2\nu+1)W_{\nu}(x)\right]
$$
and
$$\mathscr{W}''_{\nu}(x)=2^{\nu}\Gamma(\nu+1)\Gamma(\nu+2)x^{-2\nu-3}\left[x^2W''_{\nu}(x)-(4\nu+2)xW'_{\nu}(x)+(2\nu+1)(2\nu+2)W_{\nu}(x)\right],
$$
from which we get
$$\Theta_3=\frac{1}{8\gamma^4_{\nu,k}}
\left[2\nu+5+\sum_{n\ge 1}\frac{4\gamma^4_{\nu,k}}{\gamma'^4_{\nu,n}-\gamma^4_{\nu,k}}\right].
$$
This complete the proof of the equation \eqref{id3}.
\end{proof}

\begin{proof}[\bf Proof of Theorem \ref{theorem5}]
Again using \eqref{product} we have
$$\frac{z\mathscr{W}'_{\nu}(z)}{\mathscr{W}_{\nu}(z)}=-\sum_{n\geq 1}\frac{4z^4}{\gamma^4_{\nu,n}-z^4}
=-4\sum_{n\geq 1}\frac{z^4/\gamma^4_{\nu,n}}{1-z^4/\gamma^4_{\nu,n}}
=-4\sum_{n\geq 1}\sum_{m\geq 1}\frac{z^{4m}}{\gamma^{4m}_{\nu,n}}
=-4\sum_{m\geq 1}\left(\sum_{n\geq 1}\frac{1}{\gamma^{4m}_{\nu,n}}\right)z^{4m},$$
which is valid for $|z|<\gamma_{\nu,1}$ and $\nu>-1$. Hence the conclusion follows.
\end{proof}

\begin{proof}[\bf Proof of Theorem \ref{theorem6}]
The infinite product representation \eqref{product} yields
$$
\left[\log \left(x^{\frac{-2\nu-1}{4}}W_{\nu}(\sqrt[4]{x})\right)\right]'=\sum_{n\geq 1}\frac{1}{x-\gamma^4_{\nu,n}}.
$$
This gives
$$
f_{\mu,\nu}(x)=\sum_{n\geq 1}\left(\frac{1}{\gamma^4_{\nu,n}-x}-\frac{1}{\gamma^4_{\mu,n}-x}\right)+\frac{1}{16}\left(\frac{1}{(\mu+1)_3}-\frac{1}{(\nu+1)_3}\right)
$$
and hence on differentiating $m$ times we get
$$
f^{(m)}_{\mu,\nu}(x)=\sum_{n\geq 1}\left(\frac{m!}{(\gamma^4_{\nu,n}-x)^{m+1}}-\frac{m!}{(\gamma^4_{\mu,n}-x)^{m+1}}\right)\geq 0,
$$
for all $m\in \mathbb{N}$, $\mu\geq\nu>-1$ and $x\in[0,\gamma^4_{\nu,1})$.
Here we used the monotonicity of zeros of cross-product of Bessel functions \cite{abp}, namely $\nu\mapsto \gamma_{\nu,n}$ is increasing on $(-1,\infty)$ for $n\in \mathbb{N}$ fixed. Therefore, for all $n,m\in \mathbb{N}$, $\mu\geq\nu>-1$ and $x\in[0,\gamma^4_{\nu,1})$, we have
$(x-\gamma^4_{\nu,n})^{m+1}\leq (x-\gamma^4_{\mu,n})^{m+1}$ and consequently the above inequality follows. Since $f_{\mu,\nu}$ is increasing on $[0,\gamma^4_{\nu,1})$ for all $\mu\geq\nu>-1$ and in view of \eqref{rf4}, $f_{\mu,\nu}(0)=0$ we obtain that $f_{\mu,\nu}(x)\geq f_{\mu,\nu}(0)=0$. Therefore $x\mapsto f_{\mu,\nu}(x)$ is absolutely monotonic on $[0,\gamma^4_{\nu,1})$ for all $\mu\geq\nu>-1$.

Now, consider
$$h_{\nu}(x)=\left[\log \left(\frac{x^{\frac{\nu}{2}+\frac{1}{4}}e^{\frac{-x}{16(\nu+1)_3}}}{W_{\nu}(\sqrt[4]{x})}\right)\right]'=-\frac{1}{16(\nu+1)_3}+\sum_{n\geq 1}\frac{1}{\gamma^4_{\nu,n}-x}.$$
Therefore by differentiating $m$ times we have
$$
h^{(m)}_{\nu}(x)=\sum_{n\geq 1}\frac{m!}{(\gamma^4_{\nu,n}-x)^m}\geq 0
$$
for all $m\in \mathbb{N}$, $\nu>-1$ and $x\in[0,\gamma^4_{\nu,1})$. Hence $h_{\nu}$ is increasing on $[0,\gamma^4_{\nu,1})$ for all $\nu>-1$ and in view of \eqref{rf4}, $h_{\nu}(0)=0$ we obtain that $h_{\nu}(x)\geq h_{\nu}(0)=0$. This proves the absolute monotonicity of $x\mapsto h_{\nu}(x)$ on $[0,\gamma^4_{\nu,1})$ for all $\nu>-1$.

Finally, by using the fact that the exponential of a function having an absolutely monotonic derivative is absolutely monotonic, we conclude that $x\mapsto g_{\mu,\nu}(x)$ and $x\mapsto q_{\nu}(x)$ are absolutely monotonic on $[0,\gamma^4_{\nu,1})$.
\end{proof}

\begin{proof}[\bf Proof of Corollary \ref{corol6.1}]
Since $x\mapsto q_{\nu}(x)$ absolutely monotonic on $[0,\gamma^4_{\nu,1})$, it is increasing. Therefore from \eqref{product} we get
$$
q_{\nu}(x)\geq q_{\nu}(0)=2^{2\nu}\Gamma(\nu+1)\Gamma(\nu+2),
$$
which implies that
$$
W(\sqrt[4]{x})\leq \frac{x^{\frac{2\nu+1}{4}}e^{-\frac{x}{16(\nu+1)_3}}}{2^{2\nu}\Gamma(\nu+1)\Gamma(\nu+2)}.
$$
Hence by changing $x$ to $x^4$ we get the required inequality.
\end{proof}

\begin{proof}[\bf Proof of Theorem \ref{theorem7}]
{\bf a \& c.} To prove the inequality \eqref{redheffer2}, it is enough to establish the following inequality
$$
\mathcal{D}_{\nu}(x)\geq \frac{1-x^2}{1+x^2}~~\mbox{for all}~~|x|\leq \frac{\delta_{\nu}}{\alpha_{\nu,1}}.
$$
Taking into account the infinite product representation \eqref{dfprod}, we have that
\begin{equation}\label{the7-eq1}
\mathcal{D}_{\nu}(x\alpha_{\nu,1})=\frac{1-x^2}{1+x^2}\left[(1+x^2)\lim_{n\rightarrow \infty}F_{\nu,n}(x)\right],
\end{equation}
where
$$
F_{\nu,n}(x)=\prod_{k=2}^{n}\left(1-\frac{x^2\alpha^2_{\nu,1}}{\alpha^2_{\nu,k}}\right).
$$
Making use of the principle of mathematical induction we show that the following inequality
\begin{equation}\label{the7-eq2}
(1+x^2)F_{\nu,n}(x)\geq 1+\frac{x^2\alpha_{\nu,1}}{\alpha_{\nu,n}}
\end{equation}
is valid for all $\nu>-1$, $n\geq 2$ and $|x|\leq \frac{\delta_{\nu}}{\alpha_{\nu,1}}$. For $n=2$, the inequality \eqref{the7-eq2} follows from the assumption in the statement of the theorem. Namely, we have
$$
(1+x^2)F_{\nu,2}(x)-\left(1+\frac{x^2\alpha_{\nu,1}}{\alpha_{\nu,2}}\right)=\frac{x^2}{\alpha^2_{\nu,2}}\left(\Psi_{\nu}(1)-\alpha^2_{\nu,1}x^2\right)\geq 0.
$$
Now, let us assume that the inequality \eqref{the7-eq2} holds for some $m\geq 2$. Therefore
\begin{eqnarray*}
(1+x^2)F_{\nu,m+1}(x)-\left(1+\frac{x^2\alpha_{\nu,1}}{\alpha_{\nu,m+1}}\right)&=&(1+x^2)F_{\nu,m}(x)\left(1-\frac{x^2\alpha^2_{\nu,1}}{\alpha^2_{\nu,m+1}}\right)-\left(1+\frac{x^2\alpha_{\nu,1}}{\alpha_{\nu,m+1}}\right)\\
&\geq &\left(1+\frac{x^2\alpha_{\nu,1}}{\alpha_{\nu,m}}\right)\left(1-\frac{x^2\alpha^2_{\nu,1}}{\alpha^2_{\nu,m+1}}\right)-\left(1+\frac{x^2\alpha_{\nu,1}}{\alpha_{\nu,m+1}}\right)\\
&=&\frac{x^2\alpha_{\nu,1}}{\alpha_{\nu,m}\alpha^2_{\nu,m+1}}\left(\Psi_{\nu}(m)-\alpha^2_{\nu,1}x^2\right)\geq 0.
\end{eqnarray*}
Hence, by the principle of mathematical induction, inequality \eqref{the7-eq2} holds for all $n\geq 2$. Now taking limit $n\rightarrow \infty$ in \eqref{the7-eq2} we get
$$
\lim_{n\rightarrow \infty}(1+x^2)F_{\nu,n}(x)\geq \lim_{n\rightarrow \infty}\left(1+\frac{x^2\alpha_{\nu,1}}{\alpha_{\nu,n}}\right)=1,
$$
which in view of \eqref{the7-eq1} gives the inequality \eqref{redheffer1}.

To prove the inequality \eqref{redheffer3}, similar to the part {\bf a}, it is enough to prove the inequality
$$
\mathscr{W}_{\nu}(x)\geq \frac{1-x^4}{1+x^4}~~\mbox{for all}~~|x|\leq \frac{\epsilon_{\nu}}{\gamma_{\nu,1}}.
$$
Now using \eqref{product}, we have
\begin{equation}\label{the7-eq3}
\mathscr{W}_{\nu}(x\gamma_{\nu,1})=\frac{1-x^4}{1+x^4}\left[(1+x^4)\lim_{n\rightarrow \infty}G_{\nu,n}(x)\right],
\end{equation}
where
$$
G_{\nu,n}(x)=\prod_{k=2}^{n}\left(1-\frac{x^4\gamma^4_{\nu,1}}{\gamma^4_{\nu,k}}\right).
$$
Using the principle of mathematical induction we show that the inequality
\begin{equation}\label{the7-eq4}
(1+x^4)G_{\nu,n}(x)\geq 1+\frac{x^4\gamma^2_{\nu,1}}{\gamma^2_{\nu,n}}
\end{equation}
holds for all $\nu>-1$, $n\geq 2$ and $|x|\leq \frac{\epsilon_{\nu}}{\gamma_{\nu,1}}$. For $n=2$, \eqref{the7-eq4} follows from the assumption of the theorem. That is, we have
$$
(1+x^4)G_{\nu,2}(x)-\left(1+\frac{x^4\gamma^2_{\nu,1}}{\gamma^2_{\nu,2}}\right)=\frac{x^4}{\gamma^4_{\nu,2}}\left(\Omega_{\nu}(1)-\gamma^4_{\nu,1}x^4\right)\geq 0.
$$
Now, let us assume that the inequality \eqref{the7-eq4} holds for some $m\geq 2$. Therefore
\begin{eqnarray*}
(1+x^4)G_{\nu,m+1}(x)-\left(1+\frac{x^4\gamma^2_{\nu,1}}{\gamma^2_{\nu,m+1}}\right)&=&(1+x^4)G_{\nu,m}(x)\left(1-\frac{x^4\gamma^4_{\nu,1}}{\gamma^4_{\nu,m+1}}\right)-\left(1+\frac{x^4\gamma^2_{\nu,1}}{\gamma^2_{\nu,m+1}}\right)\\
&\geq &\left(1+\frac{x^4\gamma^2_{\nu,1}}{\gamma^2_{\nu,m}}\right)\left(1-\frac{x^4\gamma^4_{\nu,1}}{\gamma^4_{\nu,m+1}}\right)-\left(1+\frac{x^4\gamma^2_{\nu,1}}{\gamma^2_{\nu,m+1}}\right)\\
&=&\frac{x^4\gamma^2_{\nu,1}}{\gamma^2_{\nu,m}\gamma^4_{\nu,m+1}}\left(\Omega_{\nu}(m)-\gamma^4_{\nu,1}x^4\right)\geq 0.
\end{eqnarray*}
Consequently, by the principle of mathematical induction, inequality \eqref{the7-eq4} holds for all $n\geq 2$. Now taking the limit $n\rightarrow \infty$ in \eqref{the7-eq4} we get
$$
\lim_{n\rightarrow \infty}(1+x^4)G_{\nu,n}(x)\geq \lim_{n\rightarrow \infty}\left(1+\frac{x^4\gamma^2_{\nu,1}}{\gamma^2_{\nu,n}}\right)=1,
$$
which in view of \eqref{the7-eq3} gives the inequality \eqref{redheffer3}.

{\bf b.} Since the functions appear in the inequality \eqref{redheffer2} are even in $x$, it is enough to prove the inequality \eqref{redheffer2} for $x\in [0,\alpha_{\nu,1})$. Let us define a function $\phi_{\nu}:[0,\alpha_{\nu,1})\rightarrow \mathbb{R}$ by
$$
\phi_{\nu}(x)=\frac{3\alpha^2_{\nu,1}}{8(\nu+1)}\log \left(\frac{\alpha^2_{\nu,1}-x^2}{\alpha^2_{\nu,1}+x^2}\right)-\log \mathcal{D}_{\nu}(x),
$$
which in view of \eqref{rf3}, \eqref{ray-ineq1} and \eqref{series1} yields
\begin{eqnarray*}
\phi'_{\nu}(x)&=&-\frac{3\alpha^2_{\nu,1}}{8(\nu+1)}\cdot\frac{4x\alpha^2_{\nu,1}}{\alpha^4_{\nu,1}-x^4}-\frac{\mathcal{D}'_{\nu}(x)}{\mathcal{D}_{\nu}(x)}\\
&=&-2x\eta_{2}(\nu)\cdot\frac{\alpha^4_{\nu,1}}{\alpha^4_{\nu,1}-x^4}+\frac{2}{x}\sum_{m\geq 1}\eta_{2m}(\nu)x^{2m}\\
&=&\frac{2\eta_{2}(\nu)}{x(\alpha^2_{\nu,1}+x^2)}\left[\frac{(\alpha^2_{\nu,1}+x^2)}{\eta_{2}(\nu)}\sum_{m\geq 1}\eta_{2m}(\nu)x^{2m}-\frac{\alpha^4_{\nu,1}x^2}{\alpha^2_{\nu,1}-x^2}\right]\\
&=&\frac{2\eta_{2}(\nu)}{x(\alpha^2_{\nu,1}+x^2)}\left[\alpha^2_{\nu,1}x^2+\frac{\alpha^2_{\nu,1}}{\eta_{2}(\nu)}\sum_{m\geq 2}\eta_{2m}(\nu)x^{2m}\right.\\
&&+\left.\frac{1}{\eta_{2}(\nu)}\sum_{m\geq 1}\eta_{2m}(\nu)x^{2m+2}-\sum_{m\geq 0}\frac{x^{2m+2}}{\alpha^{2m-2}_{\nu,1}}\right]\\
&=&\frac{2\eta_{2}(\nu)}{x(\alpha^2_{\nu,1}+x^2)}\left[\frac{\alpha^2_{\nu,1}}{\eta_{2}(\nu)}\sum_{m\geq 2}\eta_{2m}(\nu)x^{2m}+\frac{1}{\eta_{2}(\nu)}\sum_{m\geq 2}\eta_{2m-2}(\nu)x^{2m}-\sum_{m\geq 2}\frac{x^{2m}}{\alpha^{2m-4}_{\nu,1}}\right]\\
&=&\frac{2\eta_{2}(\nu)}{x(\alpha^2_{\nu,1}+x^2)}\left[\sum_{m\geq 2}\left(\frac{\alpha^2_{\nu,1}\eta_{2m}(\nu)}{\eta_{2}(\nu)}+\frac{\eta_{2m-2}(\nu)}{\eta_{2}(\nu)}-\frac{1}{\alpha^{2m-4}_{\nu,1}}\right)x^{2m}\right]\\
&\geq &\frac{2\eta_{2}(\nu)}{x(\alpha^2_{\nu,1}+x^2)}\left[\sum_{m\geq 2}\left(\frac{2\alpha^2_{\nu,1}\eta_{2m}(\nu)}{\eta_{2}(\nu)}-\frac{1}{\alpha^{2m-4}_{\nu,1}}\right)x^{2m}\right]\\
&=&\frac{2\eta_{2}(\nu)}{x(\alpha^2_{\nu,1}+x^2)}\left[\sum_{m\geq 2}\left(\frac{2\alpha^{-2}_{\nu,1}\alpha^{2m}_{\nu,1}\eta_{2m}(\nu)-\eta_{2}(\nu)}{\eta_{2}(\nu)\alpha^{2m-4}_{\nu,1}}\right)x^{2m}\right]\\
&\geq &\frac{2\eta_{2}(\nu)}{x(\alpha^2_{\nu,1}+x^2)}\left[\sum_{m\geq 2}\left(\frac{2\alpha^{-2}_{\nu,1}-\eta_{2}(\nu)}{\eta_{2}(\nu)\alpha^{2m-4}_{\nu,1}}\right)x^{2m}\right]\\
&=&\frac{2\eta_{2}(\nu)}{x(\alpha^2_{\nu,1}+x^2)}\left[\sum_{m\geq 2}\left(\frac{\frac{8(\nu+1)}{3}\alpha^{-2}_{\nu,1}-1}{\alpha^{2m-4}_{\nu,1}}\right)x^{2m}\right]\\
&>&\frac{2\eta_{2}(\nu)}{x(\alpha^2_{\nu,1}+x^2)}\left[\sum_{m\geq 2}\left(\frac{8-\nu}{9(\nu+2)\alpha^{2m-4}_{\nu,1}}\right)x^{2m}\right]>0.
\end{eqnarray*}
Here in last inequality we have used the upper bound for the smallest positive zero of the Dini function (see \cite[p. 11]{ismail} with $\alpha+\nu=1$)
$$
\alpha^2_{\nu,1}=x^2_1<\frac{4(\alpha+\nu+2)(\nu+\alpha)(\nu+1)(\nu+2)}{(\alpha+\nu)^2+4\alpha+8\nu+8}=\frac{12(\nu+1)(\nu+2)}{13+4\nu}.
$$
Therefore for $\nu\in (-1,8)$, the function $\phi_{\nu}$ is increasing on $[0,\alpha_{\nu,1})$ and hence $\phi_{\nu}(x)\geq \phi_{\nu}(0)=0$ and consequently the inequality \eqref{redheffer2} holds.

Now, by using the L'Hospital rule, \eqref{rf3} and \eqref{logderivative1}, we have the limit
$$
\lim_{x\rightarrow 0}\frac{\log \mathcal{D}_{\nu}(x)}{\log \left(\frac{\alpha^2_{\nu,1}-x^2}{\alpha^2_{\nu,1}+x^2}\right)}=\lim_{x\rightarrow 0}\frac{\mathcal{D}'_{\nu}(x)}{\mathcal{D}_{\nu}(x)}\cdot\frac{x^4-\alpha^4_{\nu,1}}{4x\alpha^2_{\nu,1}}=\frac{3\alpha^2_{\nu,1}}{8(\nu+1)}=m_{\nu}.
$$
This implies that indeed the constant $m_{\nu}$ is best possible.

{\bf d.} Similar to the proof of part {\bf b} of this theorem, it is enough to prove the inequality \eqref{redheffer2'} for $x\in [0,\gamma_{\nu,1})$. Let us define a function $\Phi_{\nu}:[0,\gamma_{\nu,1})\rightarrow \mathbb{R}$ by
$$
\Phi_{\nu}(x)=\frac{\gamma^4_{\nu,1}}{32(\nu+1)_3}\log \left(\frac{\gamma^4_{\nu,1}-x^4}{\gamma^4_{\nu,1}+x^4}\right)-\log \mathscr{W}_{\nu}(x),
$$
which on using \eqref{rf4}, \eqref{ray-ineq2} and \eqref{series2} gives
\begin{eqnarray*}
\Phi'_{\nu}(x)&=&-\frac{\gamma^4_{\nu,1}}{32(\nu+1)_3}\cdot\frac{8x^3\gamma^4_{\nu,1}}{\gamma^8_{\nu,1}-x^8}-\frac{\mathscr{W}'_{\nu}(x)}{\mathscr{W}_{\nu}(x)}\\
&=&-4x^3\zeta_{4}(\nu)\cdot\frac{\gamma^8_{\nu,1}}{\gamma^8_{\nu,1}-x^8}+\frac{4}{x}\sum_{m\geq 1}\zeta_{4m}(\nu)x^{4m}\\
&=&\frac{4\zeta_{4}(\nu)}{x(\gamma^4_{\nu,1}+x^4)}\left[\frac{(\gamma^4_{\nu,1}+x^4)}{\zeta_{4}(\nu)}\sum_{m\geq 1}\zeta_{4m}(\nu)x^{4m}-\frac{\gamma^8_{\nu,1}x^4}{\gamma^4_{\nu,1}-x^4}\right]\\
&=&\frac{4\zeta_{4}(\nu)}{x(\gamma^4_{\nu,1}+x^4)}\left[\gamma^4_{\nu,1}x^4+\frac{\gamma^4_{\nu,1}}{\zeta_{4}(\nu)}\sum_{m\geq 2}\zeta_{4m}(\nu)x^{4m}\right.\\
&&+\left.\frac{1}{\zeta_{4}(\nu)}\sum_{m\geq 1}\zeta_{4m}(\nu)x^{4m+4}-\sum_{m\geq 0}\frac{x^{4m+4}}{\gamma^{4m-4}_{\nu,1}}\right]\\
&=&\frac{4\zeta_{4}(\nu)}{x(\gamma^4_{\nu,1}+x^4)}\left[\frac{\gamma^4_{\nu,1}}{\zeta_{4}(\nu)}\sum_{m\geq 2}\zeta_{4m}(\nu)x^{4m}+\frac{1}{\zeta_{4}(\nu)}\sum_{m\geq 2}\zeta_{4m-4}(\nu)x^{4m}-\sum_{m\geq 2}\frac{x^{4m}}{\gamma^{4m-8}_{\nu,1}}\right]\\
&=&\frac{4\zeta_{2}(\nu)}{x(\gamma^4_{\nu,1}+x^4)}\left[\sum_{m\geq 2}\left(\frac{\gamma^4_{\nu,1}\zeta_{4m}(\nu)}{\zeta_{4}(\nu)}+\frac{\zeta_{4m-4}(\nu)}{\zeta_{4}(\nu)}-\frac{1}{\gamma^{4m-8}_{\nu,1}}\right)x^{4m}\right]\\
&\geq &\frac{4\zeta_{4}(\nu)}{x(\gamma^4_{\nu,1}+x^4)}\left[\sum_{m\geq 2}\left(\frac{2\gamma^4_{\nu,1}\zeta_{4m}(\nu)}{\zeta_{4}(\nu)}-\frac{1}{\gamma^{4m-8}_{\nu,1}}\right)x^{4m}\right]\\
&=&\frac{4\zeta_{4}(\nu)}{x(\gamma^4_{\nu,1}+x^4)}\left[\sum_{m\geq 2}\left(\frac{2\gamma^{-4}_{\nu,1}\gamma^{4m}_{\nu,1}\zeta_{4m}(\nu)-\zeta_{4}(\nu)}{\zeta_{4}(\nu)\gamma^{4m-8}_{\nu,1}}\right)x^{4m}\right]\\
&\geq &\frac{4\zeta_{4}(\nu)}{x(\gamma^4_{\nu,1}+x^4)}\left[\sum_{m\geq 2}\left(\frac{2\gamma^{-4}_{\nu,1}-\zeta_{4}(\nu)}{\zeta_{4}(\nu)\gamma^{4m-8}_{\nu,1}}\right)x^{4m}\right]\\
&=&\frac{4\zeta_{4}(\nu)}{x(\gamma^4_{\nu,1}+x^4)}\left[\sum_{m\geq 2}\left(\frac{32(\nu+1)_3\gamma^{-4}_{\nu,1}-1}{\gamma^{4m-8}_{\nu,1}}\right)x^{4m}\right]\\
&>&\frac{4\zeta_{4}(\nu)}{x(\gamma^4_{\nu,1}+x^4)}\left[\sum_{m\geq 2}\left(\frac{\nu-\nu^2+14}{(\nu+4)(\nu+5)\gamma^{4m-8}_{\nu,1}}\right)x^{4m}\right]>0,
\end{eqnarray*}
where the last inequality follows by using the upper bound given in \eqref{ub1}. Therefore for $\nu\in (-1,r)$, the function $\Phi_{\nu}$ is increasing on $[0,\gamma_{\nu,1})$. This implies that $\Phi_{\nu}(x)\geq \Phi_{\nu}(0)=0$ and hence the inequality \eqref{redheffer2'} holds.

Now using the L'Hospital rule, \eqref{rf4} and \eqref{logderivative}, we have the limit
$$
\lim_{x\rightarrow 0}\frac{\log \mathscr{W}_{\nu}(x)}{\log \left(\frac{\gamma^4_{\nu,1}-x^4}{\gamma^4_{\nu,1}+x^4}\right)}=\lim_{x\rightarrow 0}\frac{\mathscr{W}'_{\nu}(x)}{\mathscr{W}_{\nu}(x)}\cdot\frac{x^8-\gamma^8_{\nu,1}}{8x^3\gamma^4_{\nu,1}}=\frac{\gamma^4_{\nu,1}}{32(\nu+1)_3}=n_{\nu}.
$$
This implies that indeed the constant $n_{\nu}$ is best possible.
\end{proof}

\begin{proof}[\bf Proof of Theorem \ref{theorem8}]
Since all the functions appear in the inequality \eqref{redheffer4} are even in $x$, it is enough to prove the inequality \eqref{redheffer4} for $x\in (0,r)$ for any given $r\in (0,\infty)$. Let us define a function $Q_{\nu}:(0,r)\rightarrow \mathbb{R}$ by
$$
Q_{\nu}(x)=\frac{\log \lambda_{\nu}(x)}{\log \left( \frac{r^2+x^2}{r^2-x^2}\right)}=\frac{f(x)}{g(x)}.
$$
Making use of the infinite product representation \eqref{mdfprod}, we obtain
$$
\frac{f'(x)}{g'(x)}=\frac{1}{2r^2}\sum_{n\geq 1}\frac{r^4-x^4}{\alpha^2_{\nu,n}+x^2}.
$$
Now, it is not difficult to verify that each term of the above series is decreasing on $(0,r)$ as a function of $x$. Thus, $x\mapsto \frac{f'(x)}{g'(x)}$ is decreasing on $(0,r)$ and consequently with the help of monotone form of L'Hospital's rule \cite[Lemma 2.2]{anderson}, we conclude that $x\mapsto Q_{\nu}(x)$ is decreasing on $(0,r)$. Moreover,
$$
\alpha=\lim_{x\rightarrow r}Q_{\nu}(x) < Q_{\nu}(x) < \lim_{x\rightarrow 0}Q_{\nu}(x)=\beta.
$$
This completes the proof of \eqref{redheffer4}.
\end{proof}

\begin{proof}[\bf Proof of Theorem \ref{theorem9}]
In order to prove the inequalities \eqref{redheffer7} and \eqref{redheffer8}, it is enough to consider the case  $x\in (0,j_{\nu,1})$ as  all the functions appear in \eqref{redheffer7} and \eqref{redheffer8} are even in $x$. Define a function $k_{\nu}:[0,j_{\nu,1})\rightarrow \mathbb{R}$ by
$$
k_{\nu}(x)=\log \mathcal{J}_{\nu}(x)-\log \left(\frac{j^2_{\nu,1}-x^2}{j^2_{\nu,1}+x^2}\right).
$$
Now we recall Kishore's formula \cite{kishore}
$$
\frac{x}{2}\frac{J_{\nu+1}(x)}{J_{\nu}(x)}=\sum_{m\geq 1}\sigma^{(2m)}_{\nu}x^{2m}
$$
which in view of the identity $\frac{\mathcal{J}'_{\nu}(x)}{\mathcal{J}_{\nu}(x)}=-\frac{J_{\nu+1}(x)}{J_{\nu}(x)}$ can be rewritten as
$$
\frac{x}{2}\frac{\mathcal{J}'_{\nu}(x)}{\mathcal{J}_{\nu}(x)}=-\sum_{m\geq 1}\sigma^{(2m)}_{\nu}x^{2m}.
$$
Therefore on using the above equation for $k_{\nu}(x)$ we have
\begin{eqnarray*}
k'_{\nu}(x)&=&\frac{\mathcal{J}'_{\nu}(x)}{\mathcal{J}_{\nu}(x)}+\frac{4xj^2_{\nu,1}}{(j^2_{\nu,1}+x^2)(j^2_{\nu,1}-x^2)}\\
&=&-\frac{2}{x}\sum_{m\geq 1}\sigma^{(2m)}_{\nu}x^{2m}+\frac{4x}{j^2_{\nu,1}+x^2}\sum_{m\geq 0}\frac{x^{2m}}{j^{2m}_{\nu,1}}\\
&=&\frac{2}{x(j^2_{\nu,1}+x^2)}\left[-(j^2_{\nu,1}+x^2)\sum_{m\geq 1}\sigma^{(2m)}_{\nu}x^{2m}+2\sum_{m\geq 0}\frac{x^{2m+2}}{j^{2m}_{\nu,1}}\right]\\
&=&\frac{2}{x(j^2_{\nu,1}+x^2)}\left[x^2\left(2-j^2_{\nu,1}\sigma^{(2)}_{\nu}\right)-\sum_{m\geq 2}j^2_{\nu,1}\sigma^{(2m)}_{\nu}x^{2m}-\sum_{m\geq 1}\sigma^{(2m)}_{\nu}x^{2m+2}+2\sum_{m\geq 1}\frac{x^{2m+2}}{j^{2m}_{\nu,1}}\right]\\
&=&\frac{2}{x(j^2_{\nu,1}+x^2)}\left[x^2\left(2-\frac{j^2_{\nu,1}}{4(\nu+1)}\right)+\sum_{m\geq 2}\left(\frac{2}{j^{2m-2}_{\nu,1}}-j^2_{\nu,1}\sigma^{(2m)}-\sigma^{(2m-2)}\right)x^{2m}\right]\\
&=&\frac{2}{x(j^2_{\nu,1}+x^2)}\left[x^2\left(\frac{8(\nu+1)-j^2_{\nu,1}}{4(\nu+1)}\right)+\sum_{m\geq 2}\left(\frac{(1-j^{2m}_{\nu,1}\sigma^{(2m)})+(1-j^{2m-2}_{\nu,1}\sigma^{(2m-2)})}{j^{2m-2}_{\nu,1}}\right)x^{2m}\right]\\
&\geq &0,
\end{eqnarray*}
where $\nu\geq \nu_0$. Here we have used Lemma \ref{lemma2} and the left-hand side of Rayleigh inequality \eqref{ray-ineq3}. Therefore the function $k_{\nu}$ is decreasing on $[0,j_{\nu,1})$ for all $\nu\geq \nu_0$. Consequently $k_{\nu}(x)\leq k_{\nu}(0)=0$ and hence the inequality \eqref{redheffer7} follows.

Now, taking into account of the inequality \eqref{redheffer7}, the  following inequality \cite[Theorem 3]{baricz2}
$$
\left[\mathcal{J}_{\nu+1}(x)\right]^{\nu+2}\geq \left[\mathcal{J}_{\nu}(x)\right]^{\nu+1},
$$
which is valid for all $\nu>-1$ and $x\in (-j_{\nu,1},j_{\nu,1})$ gives
$$
\frac{\mathcal{J}_{\nu+1}(x)}{\mathcal{J}_{\nu}(x)}\geq \left[\mathcal{J}_{\nu}(x)\right]^{\frac{\nu+1}{\nu+2}-1}=\frac{1}{\left[\mathcal{J}_{\nu}(x)\right]^{1/(\nu+2)}}\geq
\left(\frac{j^2_{\nu,1}+x^2}{j^2_{\nu,1}-x^2}\right)^{\frac{1}{\nu+2}}.
$$
Hence the inequality \eqref{redheffer8} is indeed true.

\end{proof}

\begin{proof}[\bf Proof of Theorem \ref{theorem10}]
Let $\nu>-1$ and $0<|x|<\alpha_{\nu,1}$. Then first we prove the following identity
\begin{equation}\label{id4}
A_{\nu}(x)=-\frac{2(\nu+1)}{3}\cdot\frac{\alpha^2_{\nu,1}-x^2}{x}\cdot\frac{\mathcal{D}'_{\nu}(x)}{\mathcal{D}_{\nu}(x)}=
\alpha^2_{\nu,1}+\frac{4(\nu+1)}{3}\sum_{m\geq 1}A_mx^{2m},
\end{equation}
where $A_m=\alpha^2_{\nu,1}\eta_{2m+2}(\nu)-\eta_{2m}(\nu)$.

To prove \eqref{id4}, we appeal to the equations \eqref{rf3} and \eqref{series1}
to obtain
\begin{eqnarray*}
A_{\nu}(x)&=&\frac{4(\nu+1)}{3}\frac{\alpha^2_{\nu,1}-x^2}{x^2}\left(-\frac{x\mathcal{D}'_{\nu}(x)}{2\mathcal{D}_{\nu}(x)}\right)\\
&=&\frac{1}{\eta_2(\nu)}\frac{\alpha^2_{\nu,1}-x^2}{x^2}\sum_{m\geq 1}\eta_{2m}(\nu)x^{2m}\\
&=&\frac{\alpha^2_{\nu,1}}{\eta_2(\nu)}\left(\eta_2(\nu)+\sum_{m\geq 2}\eta_{2m}(\nu)x^{2m-2}\right)-\frac{1}{\eta_2(\nu)}\sum_{m\geq 1}\eta_{2m}(\nu)x^{2m}\\
&=&\alpha^2_{\nu,1}+\frac{1}{\eta_2(\nu)}\sum_{m\geq 1}\left(\alpha^2_{\nu,1}\eta_{2m+2}(\nu)-\eta_{2m}(\nu)\right)x^{2m}\\
&=&\alpha^2_{\nu,1}+\frac{4(\nu+1)}{3}\sum_{m\geq 1}A_mx^{2m}.
\end{eqnarray*}
Now, for a given $N\in\mathbb{N}$, let us consider
$$
\mathscr{A}(x)=\frac{1}{x^{2n+2}}\left(\frac{A_{\nu}(x)-\alpha^2_{\nu,1}}{\frac{4(\nu+1)}{3}}-\sum_{m=1}^{n}A_mx^{2m}\right),
$$
which in view of \eqref{id4} can be rewritten as
$$
\mathscr{A}(x)=\sum_{m\geq n+1}A_mx^{2m-2n-2}=\sum_{m\geq0}A_{n+1+m}x^{2m}.
$$
Taking into account the right-hand side of \eqref{ray-ineq1}, $A_n<0$ for all $n\in \mathbb{N}$ and consequently from the above expression, $x\mapsto\mathscr{A}(x)$ is strictly decreasing on $(0,\alpha_{\nu,1})$, which implies that
$$
a=\lim_{x\rightarrow \alpha^-_{\nu,1}}\mathscr{A}(x)<\mathscr{A}(x)<\lim_{x\rightarrow 0^+_{\nu,1}}\mathscr{A}(x)=b,
$$
where $b=A_{n+1}$ and in view of the limit $
\lim_{x\rightarrow \alpha^-_{\nu,1}}A_{\nu}(x)=\frac{4(\nu+1)}{3}
$, we have
$$a=\frac{1}{\alpha^{2n+2}_{\nu,1}}\left(1-\frac{3\alpha^2_{\nu,1}}{4(\nu+1)}
-\sum_{m=1}^nA_m\alpha^{2m}_{\nu,1}\right)
$$
This completes the proof.
\end{proof}

\begin{proof}[\bf Proof of Theorem \ref{theorem11}]
Let $\nu>-1$ and $0<|x|<\gamma_{\nu,1}$. Then we need to prove the identity
\begin{equation}\label{id5}
B_{\nu}(x)=-4(\nu+1)_3\cdot\frac{\gamma^4_{\nu,1}-x^4}{x^3}\cdot\frac{\mathscr{W}'_{\nu}(x)}{\mathscr{W}_{\nu}(x)}=\gamma^4_{\nu,1}+16(\nu+1)_3\sum_{m\geq 1}B_mx^{4m},
\end{equation}
where $B_m=\gamma^4_{\nu,1}\zeta_{4m+4}(\nu)-\zeta_{4m}(\nu)$.
In order to prove \eqref{id5}, we use the equations \eqref{rf4} and \eqref{series2} and obtain
\begin{eqnarray*}
B_{\nu}(x)&=&16(\nu+1)_3\cdot\frac{\gamma^4_{\nu,1}-x^4}{x^4}\left(-\frac{x\mathscr{W}'_{\nu}(x)}{4\mathscr{W}_{\nu}(x)}\right)\\
&=&\frac{1}{\zeta_4(\nu)}\frac{\gamma^4_{\nu,1}-x^4}{x^4}\sum_{m\geq 1}\zeta_{4m}(\nu)x^{4m}\\
&=&\frac{\gamma^4_{\nu,1}}{\zeta_4(\nu)}\left(\zeta_4(\nu)+\sum_{m\geq 2}\zeta_{4m}(\nu)x^{4m-4}\right)-\frac{1}{\zeta_4(\nu)}\sum_{m\geq 1}\zeta_{4m}(\nu)x^{4m}\\
&=&\gamma^4_{\nu,1}+\frac{1}{\zeta_4(\nu)}\sum_{m\geq 1}\left(\gamma^4_{\nu,1}\zeta_{4m+4}(\nu)-\zeta_{4m}(\nu)\right)x^{4m}\\
&=&\gamma^4_{\nu,1}+16(\nu+1)_3\sum_{m\geq 1}B_mx^{4m}.
\end{eqnarray*}
Now, for a given $N\in\mathbb{N}$, consider
$$
\mathscr{B}(x)=\frac{1}{x^{4n+4}}\left(\frac{B_{\nu}(x)-\gamma^4_{\nu,1}}{16(\nu+1)_3}-\sum_{m=1}^{n}B_mx^{4m}\right),
$$
which in view of \eqref{id5} can be rewritten as
$$
\mathscr{B}(x)=\sum_{m\geq n+1}A_mx^{4m-4n-4}=\sum_{m\geq0}A_{n+1+m}x^{4m}.
$$
Using the right-hand side of \eqref{ray-ineq2}, $B_n<0$ for all $n\in \mathbb{N}$ and hence from the above expression, $x\mapsto\mathscr{B}(x)$ is strictly decreasing on $(0,\gamma_{\nu,1})$. From this we obtain
$$
r=\lim_{x\rightarrow \gamma^-_{\nu,1}}\mathscr{B}(x)<\mathscr{B}(x)<\lim_{x\rightarrow 0^+_{\nu,1}}\mathscr{B}(x)=s,
$$
where $s=B_{n+1}$ and by taking into account of the limit $\lim_{x\rightarrow \gamma^-_{\nu,1}}B_{\nu}(x)=16(\nu+1)_3$, one has $$r=\frac{1}{\gamma^{4n+4}_{\nu,1}}\left(1-\frac{\gamma^4_{\nu,1}}{16(\nu+1)_3}
-\sum_{m=1}^nB_m\gamma^{4m}_{\nu,1}\right).
$$
This completes the proof.
\end{proof}

\end{document}